\documentclass[11pt,letterpaper]{article}  
\usepackage{epsfig}
\usepackage{amsfonts}
\usepackage{amssymb}
\usepackage{amsmath}
\usepackage{float}
\newcommand{\ba}{\begin{array}}
\newcommand{\ea}{\end{array}}
\newcommand{\bal}{\begin{align}}
\newcommand{\eal}{\end{align}}
\newcommand{\be}{\begin{equation}}
\newcommand{\beqn}{\begin{equation}}
\newcommand{\eeqn}{\end{equation}}
\newcommand{\bea}{\begin{eqnarray}}
\newcommand{\eea}{\end{eqnarray}}
\newcommand{\benum}{\begin{enumerate}}
\newcommand{\eenum}{\end{enumerate}}
\newcommand{\bi}{\begin{itemize}}
\newcommand{\ei}{\end{itemize}}
\newcommand{\bean}{\begin{eqnarray*}} 
\newcommand{\eean}{\end{eqnarray*}}   

\def\qed{\hfill$\Box$\vskip0.3eM}

\newtheorem{theorem}{Theorem}
\newtheorem{defn}[theorem]{Definition}
\newtheorem{lemma}[theorem]{Lemma}
\newtheorem{corollary}[theorem]{Corollary}
\newtheorem{proposition}[theorem]{Proposition}
\newtheorem{example}[theorem]{Example}
\newtheorem{remark}{Remark}

\def\proof{{\it Proof.\ }}

\newcommand{\defeq}{\mathrel{\mathop:}=}
\newcommand{\eqdef}{=\mathrel{\mathop:}}

\begin{document}

\title{\LARGE\bf  Robust Mean Field Linear-Quadratic-Gaussian Games with Unknown $L^2$-Disturbance\thanks{A compressed version of this paper without detailed proofs has been presented at the 2015 IEEE CDC.}}

\date{}

\author{Jianhui Huang\thanks{J. Huang is with the Department of Applied Mathematics,
The Hong Kong Polytechnic University, Hong Kong
(majhuang@polyu.edu.hk). This author was supported by RGC Early Career Scheme (ECS) grant 502412P, 500613P.}
\quad and \quad Minyi Huang\thanks{M. Huang is with the School of Mathematics and Statistics, Carleton University, Ottawa, ON K1S 5B6, Canada
(mhuang@math.carleton.ca). This author was supported in part by Natural Sciences and Engineering Research Council (NSERC) of Canada under a Discovery Grant and a Discovery Accelerator Supplements Program. Part of this author's work was conducted  at the Department of Applied Mathematics, The Hong Kong Polytechnic University during March-May, 2014. Please address all correspondence to this author.
  }}

\maketitle

\begin{abstract}
This paper considers a class of  mean field  linear-quadratic-Gaussian (LQG) games with model uncertainty.  The drift term  in the dynamics of the agents contains a common unknown function. We take a robust optimization approach where a representative agent in the limiting model views  the drift uncertainty as an adversarial player. By including  the mean field dynamics in an augmented state space,  we solve two optimal control problems  sequentially, which combined with  consistent mean field approximations provides a solution to the robust game. A set of decentralized control strategies is derived by use of forward-backward stochastic differential equations (FBSDE) and  shown to be a  robust $\varepsilon$-Nash equilibrium.
\end{abstract}

\section{Introduction}

Mean field game theory provides an effective methodology for
the analysis and  strategy
design in a large population of players which are individually insignificant but collectively have  strong impact
(see e.g. \cite{HCM03,HCM07,HMC06,LL07}). A typical modeling analyzes a system of $N$ players with mean field coupling in their dynamics or costs, or both. The linear-quadratic-Gaussian (LQG)
  framework is of  particular interest since
it allows an explicit solution procedure.
Consider a large population  of $N$ agents.
The dynamics of agent $i$ are given by the stochastic differential equation (SDE)
\begin{align}\label{LP}
dx_i(t)=(A x_i(t)+B u_i(t)+G x^{(N)}(t))dt+DdW_{i}(t),\qquad t\ge 0,
\end{align}
where $x^{(N)}=({1}/{N})\sum_{i=1}^{N}x_i$ denotes the mean field coupling term.
The cost of agent $i$ is given by
\begin{eqnarray}\label{cost}
\mathcal{J}_{i}(u_i, \ldots, u_N)
=\mathbb{E}\Big[\int_0^{T}\left(|x_i-\Gamma x^{(N)}-\eta|_Q^2+u_i^TRu_i\right) dt+x_i^T(T)Hx_i(T)\Big],
\end{eqnarray}
where  we denote $|z|_Q= (z^TQz)^\frac12  $  and the symmetric matrices $Q\ge 0, H\ge 0$ and $R>0$.
 The LQG modeling framework was first developed in \cite{HCM03,HCM07} to obtain a set
 of strategies $(\hat u_1, \ldots, \hat u_N)$ such that each $\hat u_i$ only uses the local sample path information of $x_i$ and some deterministic functions reflecting the collective behavior of the agents and such that $(\hat u_1, \ldots, \hat u_N)$ is an $\varepsilon$-Nash equilibrium. There has existed a substantial body of literature adopting the LQG framework
\cite{B12,B11,KC13,LZ08,NCMH12,TZB14}.

For further literature, the reader is referred to \cite{CD13,CDL13,CL13,F14,HMC06,KLY11} for nonlinear diffusion based games and the associated SDE analysis, \cite{C12,LL07} for study of the coupled system of Hamilton-Jacobi-Bellman (HJB) and Fokker-Planck equations,  \cite{BFY13,BLP14,H10,NH12,NC13} for models containing a major player,  \cite{BSY13,DH15} for time consistent strategies in mean field games,  \cite{YMMS12} for mean field oscillator games, \cite{WZ13} for Markovian switching mean field games, \cite{CS14} for application to Bertrand and
Cournot equilibrium models, and \cite{AJW13} for a related solution notion called stationary equilibrium where players optimize assuming a steady-state long-run average for the empirical distribution of others' states. For an overview on mean field game theory, see \cite{BFY13,C14,GS13}.


  Within the traditional research on games, there has existed a fair amount of literature on model uncertainty.
For an $N$ player static game with finite action spaces and an uncertain payoff matrix,  a robust-optimization equilibrium is introduced in \cite{AB06} where each player optimizes its worst case payoff with respect to the uncertain set. A similar  method is applied
to hierarchical static games \cite{HF13}.
 Robustness has been addressed in dynamic games as well.
A linear-quadratic (LQ) game with system parameter uncertainties is presented in \cite{JP06}, and the deviation from the Nash equilibrium is estimated for a set of nominal strategies.
 Robust Nash equilibria are analyzed in \cite{BES03} for an LQ game with an unknown time-varying disturbance signal as an adversarial player. In the first case, a soft-constrained game is solved where
 the cost includes a quadratic penalty term for the disturbance. The second case introduces a hard constraint by specifying an $L^2$ bound on the disturbance function.
 The work \cite{KOH11} deals with stochastic games where the payoff and state transition probabilities contain uncertainty. The solution is developed by letting each player solve a robust Markov decision problem to optimize its worst case cost while other players' strategies are fixed.

This paper aims to address model uncertainty in the mean field LQG game context.
Specifically,  we focus on drift uncertainty by adding  to \eqref{LP} a common unknown $L^2$-disturbance $f$. A practical motivation is that in many decision problems, a large number of agents can share a common uncertainty source fluctuating with time, and examples include taxation, subsidy,
interest rates, and so on.
A direct consequence of our modeling is that this disturbance has global influence on the population. To address robustness, each agent locally views the disturbance as an adversarial player, and for this purpose we incorporate into \eqref{cost} an effort penalty term for the disturbance which in turn maximizes the resulting cost first. The agent minimizes subsequently. The framework of letting the disturbance maximize while its effort is penalized is called
the soft-constraint approach \cite{BB95,E06,BES03}. It has the advantage of analytical tractability. When a hard constraint is considered, the robust  mean field game is more difficult to tackle; see some preliminary analysis in \cite{HH13}. Regarding robustness in mean field games, a related work  is \cite{TBT13} where each agent is paired with its local disturbance as an adversarial player. The resulting solution is to replace the usual HJB equation by a Hamilton-Jacobi-Isaacs (HJI) equation in the solution.

 To design the individual strategies it is necessary to build the dynamics of the mean field (i.e. state average of the agents) evolving under the disturbance. This technique shares its spirit with the state  augmentation method in major player models \cite{H10,NH12,NC13}.
The subsequent robust optimization problem, as a minimax control problem,  leads to two optimal control problems with indefinite state weights \cite{W71}. They are different from the well known stochastic control problems with indefinite control weights \cite{CLZ98,LZ99}. We will follow a convex optimization approach to solve the two control problems via variational analysis and forward-backward stochastic differential equations (FBSDE) \cite{HP95,MY99,PW99}. Both the information structure and the solution procedure for our model are  different from \cite{TBT13} where each player and its local disturbance have access to its state and so dynamic programming is applicable.
Our main contributions are summarized as follows:

\begin{itemize}

\item We formulate a class of  mean field LQG  games where the players face a common uncertainty source, and introduce the robust optimization approach to solve two convex optimal control problems.

\item
 Decentralized strategies are obtained for the robust mean field game via a set of FBSDE.


\item The performance of the decentralized strategies  for the $N$ players is characterized as a robust $\varepsilon$-Nash equilibrium.

\end{itemize}

The rest of this paper is organized as follows.
Section \ref{sec:mod} introduces the mean field LQG game with a common disturbance and defines  the worst case cost for a player. Section \ref{sec:robust} studies the limiting robust optimization problem which leads to two optimal control problems solved sequentially by the disturbance and the representative player. The solution equation system of the mean field game is obtained in Section \ref{sec:cons} based on consistent mean field approximations. A key error estimate of the mean field approximation is developed in Section \ref{sec:error}. Section \ref{sec:nash} characterizes the set of decentralized strategies as a robust $\varepsilon$-Nash equilibrium.
 An extension of the analysis to players with random initial states is presented in Section \ref{sec:gen},
 and Section \ref{sec:con} concludes the paper.

\section{Mean Field LQG  Games with Drift Uncertainty}
\label{sec:mod}

Consider a finite time horizon $[0,T]$ for $T>0$. Suppose that $(\Omega,
\mathcal F, \{\mathcal F_t\}_{0\leq t\leq T}, \mathbb{P})$ is a complete
filtered probability space.
  Throughout this paper, we denote by $\mathbb{R}^{k}$ the $k$-dimensional Euclidean space, $\mathbb{R}^{n \times k}$ the set of all $n \times k$
matrices.
 We use   $|\cdot|$ to denote the norm of a  Euclidean space, or the Frobenius norm for matrices. For a  vector or matrix $M$, $M^T$ denotes its transpose. Let $L^{2}_{\mathcal{F}}(0, T; \mathbb{R}^{k})$ denote the space of all $\mathbb{R}^{k}$-valued $\mathcal{F}_t$-progressively measurable processes $x(\cdot)$ satisfying $\mathbb{E}\int_0^{T}|x(t)|^{2}dt<\infty$; $C([0,T];\mathbb{R}^k)$ (resp., $C^1([0,T];\mathbb{R}^k)$) is the space of all $\mathbb{R}^k$-valued functions $h(\cdot)$ defined on $[0,T]$ which are continuous (resp., continuously differentiable);   $L^{2}(0, T;   \mathbb{R}^{k})$ is the space of all $\mathbb{R}^{k}$-valued  measurable functions $h(\cdot)$
on $[0,T]$   satisfying $\int_0^{T}|h(t)|^{2}dt<\infty$, and we denote the norm $\|h\|_{L^2}= (\int_0^{T}|h(t)|^{2}dt)^{1/2}$.
Throughout the paper, we use $C$ (or $ C_1, C_2, \ldots$) to denote a generic constant which does not depend on the population size $N$ and  may vary from place to place.

\subsection{The game with a finite population}

Consider  $N$ agents (or players) denoted by $\mathcal{A}_{i}$, $1 \leq i \leq N$, respectively. The state $x_i$ of  ${\cal A}_i$ is $\mathbb{R}^n$-valued and satisfies the linear SDE
\begin{equation}
dx_i(t)=(Ax_i(t)+Bu_i(t)+Gx^{(N)}(t) +f(t)) dt +D dW_i(t), \quad 1\leq i \leq N, \label{x1}
\end{equation}
where $x^{(N)}=({1}/{N})\sum_{j=1}^{N}x_{j}$. The control $u_i$ takes its value in $\mathbb{R}^{n_1}$.
 The $\mathbb{R}^{n_2}$-valued standard  Brownian motions $\{W_i(t), 1\le i\le N\}$ are independent.
 The initial states $\{x_i(0), 1\leq i\leq N\}$ are deterministic and   their empirical mean has the limit $\lim_{N\to \infty}({1}/{N})\sum_{i=1}^N x_i(0)=m_0.$ We take $\{{\cal F}_t\}_{0\le t\le T}$ as the natural filtration  generated by the $Nn_2$-dimensional Brownian motion $(W_1(t), \ldots, W_N(t))$, and ${\cal F}= {\cal F}_T$.   The admissible control set $\mathcal{U}$ of ${\cal A}_i$ is
$$\mathcal{U}\defeq\left\{u_i(\cdot): u_i \in L^{2}_{\mathcal{F}}(0, T; \mathbb{R}^{n_1})\right\}.
$$
  Denote $u=(u_1, \ldots , u_{N})$ and  $u_{-i}=(u_1, \ldots, u_{i-1},$ $u_{i+1}, \ldots, u_{N})$.

The function $f\in L^2(0,T; \mathbb{R}^n)$ is an unknown disturbance  to characterize the model uncertainty, and  represents an influence from the common  environment for decision-making.  A natural motivation for considering deterministic disturbance is the following. Although each player ${\cal A}_i$ regards the disturbance as adversarial, it  {\it should not} be excessively pessimistic by assuming that  the latter will use the sample path information of $W_i$ to play against it, and instead only considers a deterministic $f$.

The cost functional of ${\cal A}_i$ is
\begin{align} \label{jcost}
{J}_i(u_i, u_{-i}, f)=\ &\mathbb{E}\left[\int_0^T \left(|x_i- (\Gamma x^{(N)}+\eta)|_Q^2
 + u_i^T R u_i -\frac{1}{\gamma} |f(t)|^2\right) dt+x_i^T(T)Hx_i(T)\right],
\end{align}
where the symmetric matrices $Q\geq 0$, $R>0$, $H\geq0$ and the constant $\gamma>0$.  We assume uniform
agents in the sense
that they share the same parameter datum $(A, B, G,  D; \Gamma, \eta, Q, R,\gamma, H).$ Also, to simplify the analysis, we consider constant parameters.


Due to the unknown function $f$,   $\mathcal{A}_i$ cannot  evaluate its cost even if all control policies $(u_1,\ldots, u_N)$ are known. To address this indeterminacy,  we  approach the game  from a robust optimization point of view where each agent takes  $f$ as an adversarial player. Here a soft-constraint \cite{BB95,BES03,E06} for the disturbance is adopted in that the term $-\frac{1}{\gamma}|f(t)|^2$ is included in \eqref{jcost} while  $f$ attempts to maximize. For given $(u_i, u_{-i})$,
define the worst case cost of ${\cal A}_i$ as
\begin{align*}
{ J}^{\rm wo}_i (u_i, u_{-i})=\sup_{f\in L^2(0,T; {\mathbb R}^n) }{ J}_i (u_i, u_{-i}, f).
\end{align*}

A set of strategies $(\hat u_1, \ldots, \hat u_N)$ is  a robust $\varepsilon$-Nash equilibrium for the $N$ players if for $\varepsilon \geq 0$,
\begin{align}
J_i^{\rm wo}(\hat u_i, \hat u_{-i})-\varepsilon \le  \inf_{u_i\in {\cal U}} J_i^{\rm wo}(u_i, \hat u_{-i})  \le J_i^{\rm wo}(\hat u_i, \hat u_{-i})
\label{jwo0}.
\end{align}
Our central objective is to design decentralized strategies based on the above solution notion.

\section{The Limiting Robust Optimization Problem}
\label{sec:robust}

We start by making an appropriate  approximation of the  coupling term $x^{(N)}$.
 Adding up the $N$ equations in \eqref{x1} and normalizing by $1/N$, we obtain
\begin{align}
dx^{(N)}=[ (A+G)x^{(N)} +  B  u^{(N)} +f]dt + D (1/N) \sum_{j=1}^N dW_j, \nonumber
\end{align}
where $u^{(N)} =({1}/{N})\sum_{j=1}^N u_j$.
Intuitively, from the point of view of ${\cal A}_i$,
$u^{(N)} $ may be approximated by a deterministic function $\bar u$.
Moreover, when $N\rightarrow \infty$, $({1}/{N})\sum_{j=1}^N dW_j$  vanishes due to the law of large numbers. In turn, a deterministic function  $m$ can be used  to approximate $x^{(N)}$. The above reasoning suggests to introduce the
limiting ordinary differential  equation (ODE)
\begin{align}
\dot{m} =(A+G) m +B\bar u +f,\quad m(0)=m_0.  \label{barxhue}
\end{align}


\subsection{The limiting model of the mean field game}

Consider the optimization problem of  a representative agent  ${\cal  A}_i$:
\begin{equation}\label{xxb0}
\begin{cases}
dx_i=( A x_i +Bu_i+G m_i +f  ) dt +DdW_i,\\
\dot{m_i}= (A  +G)m_i +B \bar u+ f  ,
\end{cases} 
\end{equation}
where the second equation is motivated from \eqref{barxhue}  and $m_i(0)=m_0$.  For the limiting model \eqref{xxb0}, $(W_i, x_i(0))$ is the same as in \eqref{x1}.  We reuse $(x_i, {\cal A}_i)$ to denote the state and the corresponding agent. This shall cause no  risk of confusion. Since $f$ will be determined as its worst case form depending on $x_i(0)$, $m_i$ is associated with the agent index $i$ so that it is ready as an appropriate notation for the subsequent closed-loop dynamics.
The cost functional is given by
\begin{align}
 \bar J_{i}(u_i,f) =\ &{\mathbb E} \int_0^T \left\{ | x_i-(\Gamma m_i +\eta)|^2_Q+u_i^T R u_i -\frac{1}{\gamma} |f(t)|^2 \right\} dt
+ {\mathbb E}x_i^T(T)Hx_i(T). \nonumber
\end{align}

We aim to find a solution pair $(\hat f, \hat u_i)$ such that
\begin{align}
\bar J_i (\hat u_i, \hat f)= \min_{u_i\in {\cal U}} \max_{f\in L^2(0,T; \mathbb{R}^n)} \bar J_i(u_i, f).  \label{jmm}
\end{align}
Finally, we need a consistency condition, i.e.,
$\frac1N \sum_{i=1}^N\hat u_i$  converges to  $ \bar u $ in some sense (this will be made precise in Section \ref{sec:cons}) and we look for  $\bar u\in C([0,T];\mathbb{R}^{n_1})$; the feasibility of doing so will be clear from our solution procedure.
The next part of our plan is to show that such strategies have the property in \eqref{jwo0} when applied in the game of $N$ agents.
In the following, we solve the optimization problem \eqref{jmm} in two steps.

\subsection{The control problem with respect to the disturbance}


Let $u_i\in {\cal U}$ and $\bar u\in C([0,T];\mathbb{R}^{n_1})$   be fixed.
   The optimal control problem is
\begin{align}
({\bf P1})\quad {\rm maximize}_{f\in L^2(0,T; {\mathbb R}^n)} \bar J_{i}(u_i,f).
\end{align}
Clearly ({ P1}) is equivalent to the following problem
\begin{align}
({\bf P1a})\quad {\rm minimize}_{f\in L^2(0,T; {\mathbb R}^n)} \bar J_{i}'(u_i,f) =\ &{\mathbb E} \int_0^T \left\{ -| x_i-(\Gamma m_i +\eta)|^2_Q +\frac{1}{\gamma} |f(t)|^2 \right\} dt\nonumber\\
 &- {\mathbb E}x_i^T(T)Hx_i(T). \nonumber
\end{align}

({P1a}) is an optimal control problem with  negative semi-definite state weights. We are interested in  the situation where  ({P1a}) is a strictly convex problem with a coercivity property. This ensures that the worse case disturbance is uniquely determined
by ${\cal A}_i$. The procedure below to identify conditions for ensuring  convexity is similar to \cite{LZ99}.


To study the convexity of $\bar J_{i}'$  in $f$, we construct a simpler auxiliary optimal control problem. Denote
\begin{align}
\widehat Q =(I-\Gamma)^TQ (I-\Gamma). \nonumber
\end{align}
Consider the dynamics
\begin{align}
\dot{z} = (A+G) z+g, \quad z(0)=0, \label{zg}
\end{align}
where $g\in L^2(0,T;\mathbb{R}^n)$.
The optimal control problem is
\begin{align*}
({\bf P1b})\quad {\rm minimize} \quad  \bar J_{i}''(g) = \int_0^T \left\{-z^T\widehat Q z +\frac{1}{\gamma} |g(t)|^2\right\} dt - z^T(T)Hz(T).
\end{align*}
For any $s\in {\mathbb R}$, we have $\bar J_{i}''(s g)= s^2 \bar J_{i}''(g)$, and so view $\bar J_{i}''$ as a quadratic functional of $g$.
\begin{defn} Let $F(g)$ be a real-valued functional of $g\in L^2(0,T; \mathbb{R}^n)$.
If  $F(g)\ge 0$ for all $g$,   $F$ is said to be positive semi-definite. If furthermore, $F(g)>0$ for all $g\ne 0$,  $F$ is said to be positive definite.
\end{defn}

\begin{lemma} \label{lemma:positive}
$\bar J_{i}'(u_i,f)$ is convex (resp., strictly convex) in $f$ if and only if $\bar J_{i}''(g)$ is positive semi-definite (resp., positive definite).
\end{lemma}

{\it Proof.} Let $(x_i, m_i)$ and $(x_i', m_i')$ be the state processes of \eqref{xxb0} corresponding to $(u_i, f)$ and $(u_i,f')$, respectively. Take any $\lambda_1\in [0, 1]$ and denote $\lambda_2= 1-\lambda_1$.
Then
\begin{align*}
&\lambda_1 \bar J'_{i} (u_i,  f) +\lambda_2 \bar J'_{i} (u_i,  f') - \bar J'_{i} (u_i,  \lambda_1 f+\lambda_2 f')  \\
 =\ & \lambda_1\lambda_2 \mathbb{E}\int_0^T\left\{ |x_i-x_i' -\Gamma(m_i-m_i')|_Q^2 +\frac1\gamma |f(t)-f'(t)|^2\right\} dt -\lambda_1\lambda_2
\mathbb{E}|x_i(T)-x_i'(T)|^2_H.
\end{align*}
Denote $g= f-f'$, and $z=x_i-x_i'$. Therefore, $z$ is deterministic and satisfies
\eqref{zg}. In addition, $m_i-m_i'=z $ for $t\in [0,T]$. Hence
$$
\lambda_1 \bar J'_{i} (u_i,  f) +\lambda_2 \bar J'_{i} (u_i,  f') - \bar J'_{i} (u_i,  \lambda_1 f+\lambda_2 f')
 =\lambda_1\lambda_2\bar J_{i}''(g)
$$
and the lemma follows.
 \qed


For our further existence analysis, we need to ensure $\bar J'_{i}(u_i, f)$
to be both strictly convex and coercive in $f$. For this purpose, we introduce the following assumption.\bigskip

({\bf H1}) There exists a small $\epsilon_0>0$ such that
$\bar J''_{i}(g)-\epsilon_0 \|g\|^2_{L^2}$ is positive semi-definite.

\bigskip

Note that (H1) is completely determined by the parameters $(\widehat Q, \gamma,\epsilon_0,H , T)$, and does not depend on $u_i$. Concerning  (H1), we have the following result. \bigskip
\begin{proposition}\label{lemma:ricj}
The following statements are equivalent:

\emph{(i)} \emph{(H1)} holds true on $[0, T].$

\emph{(ii)} The Riccati equation \begin{align}
\dot{P} +(A+G)^TP +P(A+G) - \gamma P^2 - \widehat Q   =0,\quad
P(T)= -H
 \label{ricP}
\end{align} has a unique solution on $[0, T].$

\emph{(iii)} For any $t \in [0, T],$
$$
\det\{[(0, I)e^{\mathcal{A}t}(0, I)^{T}]\}>0,
$$
where
$\mathcal{A}=\left(\begin{array}{cc} A+G+\gamma H & -\gamma I \\
\breve{Q} & -(A+G+\gamma H)^{T}\\ \end{array}\right)$
and $\breve{Q}=\gamma H^{2}+\widehat{Q}+(A+G)^{T}H+H(A+G).$
\end{proposition}

\proof In fact, (H1) is the uniform convexity condition proposed in \cite{sly}, and the equivalence between (i) and (ii) is a corollary of Theorem 4.6 of \cite{sly}. Moreover, (iii) $\Longrightarrow$ (ii) is given in Theorem 4.3 of \cite{MY99}. On the other hand, (ii) $\Longrightarrow$ (iii) is implied by Theorems 2.7 and 2.9 of \cite{YZ}.  \qed
\bigskip

For illustration of condition (ii), we give the following example.\begin{example} \label{example:A}
Consider
system \eqref{x1}-\eqref{jcost} with parameters $A=0.5$, $B=1$, $G=0.25$,  $Q=1$, $\Gamma=0.8$, $R=1.5$, $H=0$, $\gamma=1$.
Denote $\widehat A=A+G$.
We solve \eqref{ricP} to obtain
\begin{align}
P(t)=\frac{-\widehat Q(e^{\alpha (t-T)}- e^{-\alpha (t-T)})}{ \lambda_2 e^{\alpha(t-T) }- \lambda_1 e^{-\alpha (t-T)} }, \label{numP}
\end{align}
where
\begin{align*}
&\lambda_1= -\widehat A+\sqrt{\widehat A^2 -\gamma \widehat Q }=-0.027158, \quad \lambda_2= -\widehat A-\sqrt{\widehat A^2 -\gamma \widehat Q }=
-1.472842,\quad \\
&\alpha= \sqrt{\widehat A^2-\gamma \widehat Q}=0.722842.
\end{align*}
If $0<T<T_{\rm max}= \frac{1}{2\alpha}\log (\lambda_2/\lambda_1)=2.752198$,
 $P(t)$ given by \eqref{numP} is well defined on $[0,T]$. By the local Lipschitz continuity property of the vector field in \eqref{ricP}, $P(t)$ is the unique solution.
\end{example}Note that \eqref{ricP} is not a standard Riccati equation since the state weight matrix  $-\widehat Q$  is not positive semi-definite.  In general, the solvability of \eqref{ricP} cannot be ensured on an arbitrary time horizon. Condition (iii)  enables us to determine the solvability of \eqref{ricP} on a given time horizon. Note that condition (iii) is equivalent to $\det\{[(0, I)e^{\mathcal{A}t}(0, I)^{T}]\} \neq 0, \forall t \in [0, T]$ by noting $\det\{[(0, I)e^{\mathcal{A}t}(0, I)^{T}]\}_{t=0}=1.$ Condition (iii) is more checkable as illustrated by the following example.

\begin{example} \label{example:A2}
Consider system \eqref{x1}-\eqref{jcost} with parameters $A=-0.5$, $G=0.25$,  $Q=1$, $\Gamma=0.8$, $H=0$, $\gamma=1$. We obtain $\mathcal{A}=\left(\begin{array}{cc} -0.25 & -1 \\
0.04  & 0.25 \\
\end{array}\right)$,
 $e^{\mathcal{A}t}=\left(\begin{array}{cc} -\frac{1}{3}e^{\frac{3}{20}t}+
 \frac{4}{3}e^{-\frac{3}{20}t} & -\frac{10}{3}e^{\frac{3}{20}t}+\frac{10}{3}
 e^{-\frac{3}{20}t} \\
\frac{2}{15}e^{\frac{3}{20}t}-\frac{2}{15}e^{-\frac{3}{20}t}  &
\frac{4}{3}e^{\frac{3}{20}t}+\frac{1}{3}e^{-\frac{3}{20}t} \\ \end{array}\right),$   and 
 \begin{align}\det\{[(0, 1)e^{\mathcal{A}t}(0, 1)^{T}]\}=\frac{4}{3}
 e^{\frac{3}{20}t}+\frac{1}{3}e^{-\frac{3}{20}t}>0, \quad \quad \forall t\ge 0.
\end{align}
Thus for any $T>0$, \eqref{ricP} admits a unique solution on $[0, T].$ Therefore, \emph{(H1)} holds true on $[0, T].$
\end{example}

\begin{lemma} \label{lemma:coer}
Assume \emph{(H1)}. Then $\bar J'_{i}(u_i, f)$  is strictly convex in $f$. Moreover, $\bar J'_{i}(u_i, f)$  is  coercive in $f$ and, in particular, there exists a constant $C_{u_i, x_i(0)}$ depending on $(u_i, x_i(0))$ such that
 $$
\bar J'_{i} (u_i, f) \ge \frac{\epsilon_0}{2} \|f\|_{L^2}^2 -C_{u_i, x_i(0)}.
$$
\end{lemma}

\proof
Since $\bar J_i''(g)-\epsilon_0 \|g\|_{L^2}^2$ is positive semi-definite by (H1), $\bar J_i''(g)$ is positive definite. By Lemma \ref{lemma:positive}, $\bar J'_i(u_i,f)$  is  strictly convex in $f$.
 Following the method in proving Lemma \ref{lemma:positive}, we can further show that $\chi(f)\defeq \bar J'_i(u_i,f )-\epsilon_0 \|f\|^2_{L^2}$ is  convex in $f$.
 By \eqref{xxb0} and direct estimates, we can show
$$
\sup_{\|f\|_{L^2}\le 1}|\chi( f) |\le C_{0,u_i, x_i(0)},
$$
where the constant $C_{0,u_i, x_i(0)}$ depends on $(u_i, x_i(0))$.
Now consider $f$ with $\|f\|_{L^2}\ge 1$. Define $f_1=\frac{f}{\|f\|_{L^2}}$. The convexity of $\chi(f)$ implies
\begin{align}
\chi(f_1) \le \frac{1}{\|f\|_{L^2}}  \chi(f)
+\frac{\|f\|_{L^2}-1}{\|f\|_{L^2}} \chi(0)\le \frac{1}{\|f\|_{L^2}}  \chi(f)
+C_{0,u_i, x_i(0)}. \label{chi}
\end{align}
Consequently, for $\|f\|_{L^2} \ge  1$, \eqref{chi} gives
$$
\chi(f) \ge -2C_{0,u_ix_i(0)} \|f\|_{L^2}.
$$
Hence for any $f$,
$
\chi(f)\ge -C_{0,u_i, x_i(0)}(2\|f\|_{L^2}+1).
$
It follows that
\begin{align*}
\bar J'_i (u_i, f)& = \chi(f)+\epsilon_0\|f \|^2_{L^2}\\
&\ge {\epsilon_0}
\|f\|^2_{L^2}    -C_{0,u_i, x_i(0)}(2\|f\|_{L^2}+1)\\
&\ge \frac{\epsilon_0}{2}
\|f\|^2_{L^2}  -C_{u_i, x_i(0)}
\end{align*}
for some constant $C_{u_i, x_i(0)}$. \qed

\begin{theorem} \label{theorem:p1}
Suppose that \emph{(H1)} holds and let $u_i\in {\cal U}$  and $\bar u$ be fixed.
Then

\emph{(i)} $\bar J'_{i}(u_i, f)$  has a unique minimizer $\hat f$, or equivalently, $\bar J_i(u_i, f)$ has a unique maximizer $\hat f$;

\emph{(ii)} there exists a unique solution $(x_i,m_i, p_i)\in L_{\cal F}^2(0,T;\mathbb{R}^n)
\times L^2(0,T;\mathbb{R}^{2n})  $  to  the equation system
 \begin{equation}\label{bp}
 \begin{cases}
dx_i =( Ax_i +Bu_i +G m_i+\gamma  p_i) dt +D dW_i,\\
 \dot {m_i} = (A + G)m_i +B\bar u +\gamma  p_i,\\
 \dot{ p_i} = -(A+G)^T  p_i   -(I-\Gamma)^T Q [\mathbb{E}x_i- (\Gamma m_i+\eta)],
\end{cases}
\end{equation}
where  $m_i (0)= m_0$ and $ p_i(T)=H\mathbb{E}x_i(T)$, and furthermore  $\hat f=\gamma  p_i$.
\end{theorem}

\proof
(i) By Lemma \ref{lemma:positive},  $\bar J'_{i}$ is strictly convex and coercive.
In addition, $\bar J'_{i}$  is continuous in $f$. Hence there exists a unique $\hat f$ such that  $\bar J'_{i}(u_i, \hat f)=\inf_f\bar J'_{i}(u_i, f)$ \cite[Chap. 7]{KZ05}, \cite{L69}.

(ii) We start by establishing existence.  Let the optimal state-control pair be denoted by $(x_i, m_i, \hat f)$, which is uniquely determined. We have the relation
\begin{align}
&dx_i =( Ax_i +Bu_i +G m_i+\gamma \hat f) dt +D dW_i, \label{xfh}\\
& \dot{m_i} = (A + G)m_i +B\bar u +\gamma \hat f, \label{xbfh}
\end{align}
where $m_i(0)=m_0$.
By using $(x_i, m_i)$, we obtain a unique solution $ p_i$ from
\begin{equation}
\dot{p_i} = -(A+G)^T  p_i   -(I-\Gamma)^T Q [\mathbb{E}x_i- (\Gamma m_i+\eta)]  , \label{pb}
\end{equation}
where $ p_i(T) =H\mathbb{E}x_i(T)$.

Now we consider another control $f=\hat f+\tilde f\in L^2(0,T;\mathbb{R}^n)$ in place of $\hat f$.
Let $\tilde x_i$ and $\tilde m_i$ be the first variations of $x_i$ and $m_i$, respectively, which result from the variation $\tilde f$ for $\hat f$.
Then we have $\tilde x_i= \tilde m_i$ for all $t\in [0,T]$  and
\begin{align*}
\frac{d \tilde x_i }{dt}= ( A  +G) \tilde x_i +\tilde f , \quad \tilde x_i(0)=0. \end{align*}

Since $\bar J_{i}'$ has  a minimum at $(x_i, m_i, \hat f)$,
the first variation of the cost satisfies
\begin{align}
0=\frac{\delta \bar J'_{i}}{2}=\mathbb{E}\int_0^T\left\{ -[ x_i-(\Gamma m_i+\eta)]^T Q (I-\Gamma)\tilde x_i   +\frac{1}{\gamma} \hat f^T \tilde f  \right\} dt  -\mathbb{E}x_i^T(T) H\tilde x_i(T). \label{djp}
\end{align}
On the other hand,
\begin{align}
\frac{d}{dt} ( p_i^T \tilde x_i) &= \tilde x_i^T \dot{p_i} +  p_i^T\frac{ d\tilde x_i}{dt} \nonumber\\
 &=- [\mathbb{E}x_i- (\Gamma m_i+\eta)]^T Q (I-\Gamma)\tilde x_i  + p_i^T \tilde f . \label{pxode}
\end{align}
Integrating both sides of \eqref{pxode} and invoking \eqref{djp}, we obtain
\begin{align}
 p_i^T(T)\tilde x_i(T)=\int_0^T \left( p_i^T  \tilde f -\frac{1}{\gamma} \hat f^T \tilde f\right)dt
+\mathbb{E}x_i^T(T) H\tilde x_i(T). \label{pxT}
\end{align}
Recalling $ p_{i}(T)= H \mathbb{E}x_i(T)$, since $\tilde f$ is arbitrary,  it follows  from  \eqref{pxT} that
$$
\hat f=\gamma  p_i
$$
for a.e. $t\in [0,T]$. Therefore, $(x_i, m_i,  p_i)$ determined by \eqref{xfh}-\eqref{pb} is a solution to \eqref{bp}.

We proceed to show uniqueness. Suppose that $(x_i', m_i',  p_i')$ is another solution of \eqref{bp}. Set the control $f'=\gamma  p_i'$.
It is straightforward to show that the first variation of $\bar J_{i}'$ at the state control  pair $(x_i', m_i', f')$ is zero. Since $\bar J_{i}'$ is strictly convex, this implies that $(x_i', m_i', f')$ is the unique optimal state-control pair and so coincides with  $(x_i, m_i , \hat f)$ where $(x_i, m_i)$ is the optimal state process determined from \eqref{xfh}-\eqref{pb}. This further implies $ p_i'= p_i$. So  uniqueness follows. The last part of (ii) is now obvious.
\qed


\subsection{The control problem of player ${\cal A}_i$}

Assume that (H1) holds. This will ensure that all the equation systems in this section have a well defined solution.
The dynamics are given by
\begin{align} \label{bp2}
\begin{cases}
dx_i =( Ax_i +Bu_i +G m_i+\gamma  p_i) dt +D dW_i,\\
 \dot{m_i }= (A + G)m_i +B\bar u +\gamma  p_i,\\
 \dot{ p_i} = -(A+G)^T  p_i   -(I-\Gamma)^T Q [\mathbb{E}x_i- (\Gamma m_i+\eta)],
\end{cases}
\end{align}
where  $m_i(0)= m_0$ and  $ p_i(T)=H{\mathbb E}x_i(T)$.
The optimal control problem is
\begin{align}
({\bf P2})\quad {\rm minimize}_{u_i\in L_{\cal F}^2(0,T; {\mathbb R}^{n_1})} \bar J_i(u_i,\hat f_{u_i}) =\ & {\mathbb E} \int_0^T \left\{ | x_i-(\Gamma m_i +\eta)|^2_Q+u_i^T R u_i -\gamma | p_i(t)|^2 \right\} dt \nonumber \\
 & + {\mathbb E}x_i^T(T)Hx_i(T).  \nonumber
\end{align}
Here we have taken $\hat f_{u_i}=\gamma  p_i$ which depends on $u_i$. We may simply write $\bar J_i(u_i)$.
This is again a linear quadratic optimal control problem with indefinite weight for the state vector $(x_i, m_i,  p_i)$. Note that a perturbation in $u_i$ will cause a change of the mean term ${\mathbb E}x_i$. So this is essentially a mean field type optimal control problem; see related work \cite{AD10,Y13}.

We continue to identify conditions under which ({P2}) is  strictly convex and coercive.
These conditions will be characterized by using an auxiliary control problem
with dynamics
\begin{equation} \label{zzq}
\begin{cases}
\dot{z_i} = A z_i+ B \nu_i +G z  +\gamma q ,\\
\dot{ z} =(A+G) z +\gamma  q , \\
 \dot{ q} = -(A+G)^T q  -(I-\Gamma)^T Q (z_i -\Gamma z) ,
\end{cases}
\end{equation}
where $z_i(0)= z(0)=0$ and $q(T)=Hz_i(T)$. The control $\nu_i\in L^2 (0,T; \mathbb{R}^{n_1})$.
The optimal control problem is
\begin{align}
{\bf (P2a)} \quad   {\rm minimize} \quad \bar J_i^a (\nu_i) =\int_0^T \left\{| z_i - \Gamma z|_Q^2 +\nu_i^TR\nu_i -\gamma | q(t)|^2\right\} dt +
 | z_i(T)|_H^2.
\end{align}
 We may view this as a deterministic optimal control problem with two point boundary value conditions for the state trajectory.
We say $\bar J_i^a$  is positive semi-definite if $\bar J_i^a(\nu_i)\ge 0$ for all $\nu_i$; if furthermore, $\bar J_i^a(\nu_i)>0$ whenever $\nu_i\ne 0$, we say $\bar J_i^a $ is positive definite.
In order to have a well defined optimal control problem, we need to show that \eqref{zzq} has a unique solution.
\begin{lemma} \label{lemma:euzzq}
Assume \emph{(H1)}. For each $\nu_i$, there exists a unique solution $(z_i, z, q)\in C^1([0,T]; \mathbb{R}^{3n})$ to \eqref{zzq}.
\end{lemma}

\proof
Indeed, by taking $u_i=0$ and $u_i=\nu_i\in L^2 (0,T; \mathbb{R}^{n_1})$ in \eqref{bp2},  we obtain two solutions $(x_i^0, m_i^0,  p_i^0)$ and $(x_i^{\nu_i},  m_i^{\nu_i},  p_i^{\nu_i})$, respectively. It is easy to show that $(z_i, z, q)\defeq(x_i^{\nu_i}-x_i^0, m_i^{\nu_i}-m_i^0,  p_i^{\nu_i}- p_i^0)$ is a solution of \eqref{zzq} by observing that $x_i^{\nu_i}-x_i^0$ is deterministic.

If  there exist two different solutions to \eqref{zzq} for some $\nu_i$,
then we can construct two different solutions to \eqref{bp2} for a given $u_i$, which is a contradiction to Theorem \ref{theorem:p1}. \qed

\begin{lemma} \label{lemma:p2cv}
$\bar J_i(u_i)$ is convex (resp., strictly convex) in $u_i\in {\cal U}$ if and only if $\bar J_i^a(\nu_i)$ is positive semi-definite (resp., positive definite).  \end{lemma}

{\it Proof}. See appendix A.  \qed
\vskip 0.3eM

We introduce the following assumption. \bigskip

{\bf (H2)} There exists a small constant $\delta_0>0$ such that $\bar J_i^a(\nu_i)-\delta_0 \|\nu_i\|^2 \ge 0$  for all $\nu_i\in L^2(0,T; \mathbb{R}^{n_1})$.

\subsection{Representation of the quadratic functional}

We intend to find an expression of $\bar J_i^a(\nu_i)$ so that (H2) can be characterized in a more explicit form.
A change of coordinates  will make the computation more convenient.
Define
$\check z= z_i-z$. Then \eqref{zzq} becomes
\begin{align} \label{czzq}
\begin{cases}
\dot{\check{z}} = A \check z+ B \nu_i ,\\
 \dot{z} =(A+G) z +\gamma  q, \\
 \dot{ q} =-\widehat Q z -(A+G)^T q -(I-\Gamma)^T Q \check z   ,
\end{cases}
\end{align}
where $\check z(0)= z(0)=0$ and $q(T)=H(\check z(T)+ z(T))$.

Define the Hamiltonian matrix
$$
{\cal H}=\left[
\begin{array}{cc}
A+G & \gamma I \\
-\widehat Q   &  -(A+G)^T
\end{array}
 \right]
$$
and the matrix ODE $\dot\Phi (t)= {\cal H} \Phi(t)$ where $\Phi(0)=I$. Denote the partition
$$
\Phi(t)= \left[
\begin{array}{cc}
\Phi_{11}(t) & \Phi_{12}(t) \\
\Phi_{21}(t)  &\Phi_{22}(t)
\end{array}
\right],
$$
where each submatrix $\Phi_{ij}$ is an $n\times n$ matrix function.

We have
\begin{align}
\check z(t)= \int_0^t e^{A(t-\tau)} B\nu_i(\tau) d\tau . \label{zchk}
\end{align}
By solving  $(z,q)$ in \eqref{czzq}, we obtain
\begin{align*}
&z(t)= \Phi_{12}(t)q(0) -\int_0^t \Phi_{12}(t-s)  (I-\Gamma)^T Q\check z(s)ds,\\
& q(t)= \Phi_{22}(t) q(0)-\int_0^t \Phi_{22}(t-s)  (I-\Gamma)^T Q\check z(s)ds,
\end{align*}
where  $q(0)$ is to be determined.
At the terminal  time,
$$
z(T)= \Phi_{12}(T)q(0) -\int_0^T\Phi_{12}(T-s)  (I-\Gamma)^T Q\check z(s)ds
$$
and
\begin{align*}
q(T) &= \Phi_{22}(T) q(0)-\int_0^T \Phi_{22}(T-s)  (I-\Gamma)^T Q\check z(s)ds\\
   & = H\check z(T)+H \Phi_{12}(T)q(0) -H\int_0^T\Phi_{12}(T-s)  (I-\Gamma)^T Q\check z(s)ds,
\end{align*}
where the second equality is due to the terminal condition of $q$.
It follows that
\begin{align}
[\Phi_{22}(T)-H \Phi_{12}(T)]  q(0)= H\check z(T)+ \int_0^T[ \Phi_{22}(T-s)-H \Phi_{12}(T-s) ] (I-\Gamma)^T Q\check z(s)ds. \label{Pq0}
\end{align}

\begin{proposition}
If \emph{(H1)} holds, $\Phi_{22}(T)- H\Phi_{12}(T)$ is nonsingular.
\end{proposition}

\proof Under (H1),  \eqref{czzq} has a unique solution by Lemma  \ref{lemma:euzzq}, and accordingly, $q(0)$ is uniquely determined. If $\Phi_{22}(T)-H\Phi_{12} (T)$ is singular, we may find two different solutions of $q(0)$ from \eqref{Pq0} which further give two different solutions to \eqref{czzq}, leading to a contradiction.  Hence, $\Phi_{22}-H\Phi_{12}(T)$ is nonsingular.  \qed

By solving $q(0)$ in \eqref{Pq0} and further eliminating $\check z$, we write
$z$ and $q$ as integrals depending on $\nu_i$.
Define the linear operator
$$
[{\cal L}(\nu_i)](t)=\left[
\begin{array}{c}
\check z(t)\\
z(t)\\
q(t)
\end{array}
 \right] .
$$
By standard estimates we can show that ${\cal L}$ is a linear and bounded operator from $L^2(0,T; \mathbb{R}^{n_1})$ to $L^2(0,T; \mathbb{R}^{3n})$. Let ${\cal L}^*$ be its adjoint operator from $ L^2(0,T; \mathbb{R}^{3n})  $ to $L^2(0,T; \mathbb{R}^{n_1}) $.
 Define the operator
$$
{\cal L}_T \nu_i = \check z(T)+z(T).
$$
It can be shown that ${\cal L}_T$ is a linear and bounded operator from $ L^2(0,T; \mathbb{R}^{n_1})  $ to $ \mathbb{R}^n$. Let ${\cal L}^*_T$ be its adjoint operator.
Now $\bar J_i^a $ may be represented in terms of the inner product on $L^2(0,T;\mathbb{R}^{n_1})$:
\begin{align}
\bar J_i^a(\nu_i) = \langle\Theta \nu_i, \nu_i\rangle+\langle R\nu_i , \nu_i\rangle+
\langle\Theta_T \nu_i, \nu_i\rangle, \label{jrep}
\end{align}
where
$$
\Theta
\nu_i=
{\cal L}^*\left[ \begin{array}{ccc}
 Q & Q(I-\Gamma) &0 \\
(I-\Gamma)^TQ  & \widehat Q & 0 \\
0&0 & -\gamma I
\end{array} \right]
 {\cal L} \nu_i, \quad \Theta_T \nu_i= {\cal L}_T^*H{\cal L}_T \nu_i.
$$

\begin{proposition} \label{lemma:llr} \

\emph{(i)} $\bar J_i(u_i)$  is  convex in $u_i\in {\cal U}$ if and only if $\langle(\Theta +\Theta_T +R)\nu_i, \nu_i\rangle\ge 0$ for all $\nu_i\in L^2(0,T; \mathbb{R}^{n_1})$.

\emph{(ii)} \emph{(H2)} holds if and only if there exists $\delta_0>0$ such that
$\langle(\Theta +\Theta_T+R)\nu_i, \nu_i\rangle\ge \delta_0 \|\nu_i\|_{L^2}^2$ for all $\nu_i\in L^2(0,T; \mathbb{R}^{n_1})$.
\end{proposition}

\proof (i) follows from Lemma \ref{lemma:p2cv} and the representation
\eqref{jrep}. (ii) follows from \eqref{jrep}. \qed

The criterion in part (ii) of Proposition \ref{lemma:llr} still involves the operators $\Theta$ and $\Theta_T$ on an infinite dimensional
space. Here we give a sufficient condition to endure (H2) based on some more computable parameters. It is clear that
$\langle(\Theta +\Theta_T+R)\nu_i, \nu_i\rangle\ge
\int_0^T (|\nu_i(t)|^2_R -\gamma |q(t)|^2) dt.$
 For simplicity, we only consider the case
 $H=0$, and simple computations  lead to
\begin{align*}
q(t) =\ &  \Phi_{22}(t)\Phi_{22}^{-1}(T) \int_0^T \Phi_{22}(T-s)
(I-\Gamma)^TQ\int_0^se^{A(s-\tau)} B\nu_i (\tau) d\tau d s  \\
&-\int_0^t \Phi_{22}(t-s)(I-\Gamma)^TQ\int_0^se^{A(s-\tau)}
B\nu_i (\tau) d\tau d s=: q_1(t)-q_2(t) .
\end{align*}

Denote
$
b_1= \sup_{0\le t\le T}|\Phi_{22}(t)|$,
$ b_2=\sup_{0\le t \le T} |\Phi_{22}(t)\Phi_{22}^{-1}(T)|$, $b_3= |Q(I-\Gamma)|$,
$b_4=\int_0^T |e^{As} B |ds$ and $b_5=\sup_{0\le t\le T}|e^{At}B|$.
By exchanging the order of integration in $q_1$ and $q_2$, it is easy to show
$$
|q_1(t)|^2\le (b_1b_2b_3b_4)^2 {T} \int_0^T \nu_i^2(s)ds, \quad
|q_2(t)|\le b_1b_3b_5 \int_0^t (t-\tau) |\nu_i(\tau)| d\tau,
$$
which further gives
\begin{align}
\int_0^T |q(t)|^2 dt\le C_q \int_0^T |\nu_i(t)|^2 dt, \label{Cqnu}
\end{align}
where
$C_q= 2(b_1b_2b_3b_4)^2 T^2
+ \frac{1}{6}(b_1b_3b_5)^2T^4$.
For the case $H=0$,  (H2) holds  whenever $R> \gamma C_qI$.

\subsection{The solution of (P2)}

\label{sec:sub:p2}

Let $\bar u\in C([0,T];\mathbb{R}^{n_1})$ be fixed.
 \begin{lemma} \label{lemma:uy}
Assume \emph{(H1)}-\emph{(H2)}.
Then
   \emph{(P2)} has a unique optimal state-control pair of the form $(x_i, m_i,  p_i, \hat u_i)$ satisfying
\begin{equation} \label{xxp}
\begin{cases}
dx_i =( Ax_i +B\hat u_i +G m_i+\gamma  p_i) dt +D dW_i,\\
 \dot{ m_i} = (A + G)m_i +B\bar u +\gamma  p_i,\\
 \dot{ p_i }= -(A+G)^T  p_i   -(I-\Gamma)^T Q [\mathbb{E}x_i- (\Gamma  m_i+\eta)],
\end{cases}
\end{equation}
where $ p_i(T)= H{\mathbb E} x_i(T)$.
Furthermore,  the backward stochastic differential equation \emph{(BSDE)}
\begin{equation}\label{byi}
\begin{cases}
dy_i=\left\{-A^T y_i +Q[x_i-(\Gamma m_i+\eta)]\right\} dt +\zeta_i dW_i,\\
y_i(T)=- H x_i(T)
\end{cases}
\end{equation}
has a unique solution  $(y_i, \zeta_i)\in L_{\cal F}^2(0,T; \mathbb{R}^{2n}) $  and
\begin{equation}
\hat u_i= R^{-1} B^T y_i. \label{uyi}
\end{equation}
\end{lemma}

\proof
Under (H2),
by adapting Lemma \ref{lemma:p2cv} to the auxiliary control problem with cost functional
$ \bar J_i(u_i)-\delta_0\mathbb{E}\int_0^T |u_i|^2dt$, we can show   that $ \bar J_i(u_i)-\delta_0\mathbb{E}\int_0^T |u_i|^2dt$  is convex in $u_i$. By the method in proving Lemma \ref{lemma:coer}, we can further show that $\bar J_i$ is strictly convex and coercive in $u_i$.
Hence
(P2)  has a unique optimal state-control pair $(x_i, m_i,  p_i, \hat u_i)$
which minimizes $\bar J_i(u_i)$.

 Given $(x_i, m_i, p_i, \hat u_i)$, \eqref{byi} is a standard linear BSDE and so has a unique solution $(y_i,\zeta_i)$. Further define
the BSDE
\begin{equation}
dy= \left\{-G^T y_i -(A+G)^T y - \Gamma^TQ[x_i-(\Gamma m_i +\eta)]\right\} dt +\zeta dW_i, \nonumber
\end{equation}
where $y(T)=0$. It also has a unique solution $ (y, \zeta)\in L^2_{\cal F} (0,T; \mathbb{R}^{2n})$.
It can be checked that
$$
\frac{d}{dt}[{\mathbb E}(y+y_i)+ p_i ]= -(A+G)^T [{\mathbb E}(y+y_i)+ p_i]
$$
and ${\mathbb E}(y(T)+y_i(T))+ p_i(T)=0$. So
\begin{align}
{\mathbb E}(y_i+y)+ p_i=0 \label{eyy}
\end{align} for all $t\in [0,T]$.

Let $\hat u_i$ be replaced by $\hat u_i+\tilde u_i\in L^2_{\cal F} (0,T; \mathbb{R}^{n_1})$ in \eqref{xxp}, and the resulting solution  be denoted by  $(x_i+\tilde x_i, m_i+\tilde m_i, p_i+ \tilde p_i)$, which exists and is unique by Theorem \ref{theorem:p1}. It follows that
\begin{equation*}
\begin{cases}
\dot{\tilde x}_i = A\tilde x_i +B\tilde u_i +G \tilde m_i+\gamma \tilde p_i,\\
 \dot{ \tilde m}_i = (A + G)\tilde m_{i}  +\gamma \tilde p_i,\\
\dot{ \tilde p}_i = -(A+G)^T \tilde p_i   -(I-\Gamma)^T Q (\mathbb{E}\tilde x_i- \Gamma \tilde m_i),
\end{cases}
\end{equation*}
where $\tilde x_i(0)=\tilde m_i(0)=0$ and $\tilde p_i(T)=H {\mathbb E} \tilde x_i(T)$.
The first variation of $\bar J_i$ about $\hat u_i$ satisfies
\begin{equation}
0=\frac{\delta \bar J_i}{2}= {\mathbb E}\int_0^T\left\{(\tilde x_i-\Gamma \tilde m_i)^TQ[ x_i-(\Gamma m_i +\eta)]+ \tilde u_i^T R \hat u_i -\gamma \tilde  p_i^T  p_i\right\}dt  +{\mathbb E} \tilde x_i^T(T) H x_i(T). \label{dj0}
\end{equation}
By applying Ito's formula to $\tilde x_i^Ty_i$, we obtain
\begin{align}
{\mathbb E}\tilde x_i^T(T) y_i(T)-{\mathbb E}\tilde x_i^T(0) y_i(0) ={\mathbb E}\int_0^T\left\{\tilde x_i^T Q [x_i-(\Gamma m_i+\eta)] + y_i^T (B\tilde u_i +G \tilde m_i +\gamma \tilde p_i) \right\}dt .\nonumber
\end{align}
Similarly,
\begin{align}
{\mathbb E}\tilde m_i^T(T) y(T)-{\mathbb E}\tilde m_i^T(0) y(0)={\mathbb  E} \int_0^T\left\{ \gamma y^T \tilde p_i - \tilde m_i^T  (G^T y_i +\Gamma^T Q[x_i-(\Gamma m_i+\eta)])\right\}dt. \nonumber
\end{align}
Therefore, adding up the two equations yields
\begin{equation} \label{exhy}
-{\mathbb E}\tilde x_i^T(T)H x_i(T)={\mathbb E}\int_0^T \left\{ (\tilde x_i -\Gamma \tilde m_i)^TQ [x_i-(\Gamma m_i+\eta)] + y_i^T B \tilde u_i+ \gamma (y+y_i)^T \tilde p_i\right\}dt.
\end{equation}
By \eqref{dj0} and \eqref{exhy},
\begin{equation}
{\mathbb E} \int_0^T[ \tilde u_i^T R \hat u_i -\gamma \tilde  p_i^T  p_i- \tilde u_i^T B^T y_i- \gamma \tilde p_i^T(y+y_i)]dt=0.\nonumber
\end{equation}
Note that by \eqref{eyy},
$$
{\mathbb E}\int_0^T \tilde p_i^T (p_i +y+y_i)dt = \int_0^T \tilde p_i^T [ p_i+ {\mathbb E}( y+y_i)]dt =0.
$$
Hence,
$$
{\mathbb E} \int_0^T \tilde u_i^T (R \hat u_i -  B^T y_i) dt =0.
$$
Since $\tilde u_i\in L_{\cal F}^2(0,T;{\mathbb R}^{n_1})$ is arbitrary, \eqref{uyi} follows.  \qed


After substituting $\hat u_i=R^{-1}B^Ty_i$ into \eqref{xxp},
we form the  equation system
\begin{equation}\label{fbsdeubar}
\begin{cases}
dx_i =( Ax_i +BR^{-1}B^Ty_i +G m_i+\gamma  p_i) dt +D dW_i,\\
 \dot{m_i} = (A + G) m_i +B\bar u +\gamma  p_i,\\
 \dot{ p_i} = -(A+G)^T  p_i   -(I-\Gamma)^T Q [\mathbb{E}x_i- (\Gamma m_i+\eta)],\\
dy_i=\left\{-A^T y_i +Q[x_i-(\Gamma m_i+\eta)]\right\} dt +\zeta_i dW_i,
 \end{cases}
\end{equation}
where $x_i(0)$ is given, $m_i(0)=m_0$, $ p_i(T)=H{\mathbb E} x_i(T)$, and $y_i(T)=-Hx_i(T)$. This equation system consists of 2 forward equations and 2 backward equations. It is clear that the solution of the optimal control problem (P2) satisfies the above FBSDE. A natural question is whether this FBSDE's solution completely determines the optimal control. This is answered by the next theorem.
Denote
$$S[0,T]=L^2_{\cal F}(0,T;\mathbb{R}^n) \times C^1([0,T]; \mathbb{R}^{2n})\times L^2_{\cal F}(0,T;\mathbb{R}^{2n}) .$$

\begin{theorem} \label{theorem:ror}
Assume \emph{(H1)}-\emph{(H2)}. Then the \emph{FBSDE} \eqref{fbsdeubar} has a unique solution  $(x_i,  m_i,  p_i, y_i, \zeta_i)\in S[0,T]$  and the optimal control for \emph{({P2})}  is given by $\hat u_i =R^{-1} B^T y_i$.
\end{theorem}

\proof
We solve (P1) first and  (P2) next to determine $\hat u_i$.
By  Lemma \ref{lemma:uy}, we obtain $(x_i,  m_i,  p_i, y_i, \zeta_i)$
to satisfy \eqref{xxp}-\eqref{byi} and
$
\hat u_i= R^{-1} B^T y_i.
$
Obviously, $(x_i, m_i,  p_i, y_i, \zeta_i)$ satisfies \eqref{fbsdeubar}.

We continue to show uniqueness. Suppose that $(x_i, m_i,  p_i, y_i, \zeta_i)$
and $(x_i', m_i',  p_i', y_i', \zeta_i')$ are two  solutions of \eqref{fbsdeubar}. Define $  \check u_i= R^{-1} B^T y_i$ and $ u_i'= R^{-1} B^T y_i'$ which are both well-determined elements in $L^2_{\cal F}(0,T;\mathbb{R}^{n_1})$.
In particular, we have
\begin{equation} \label{x}
\begin{cases}
dx_i =( Ax_i +B\check u_i +G m_i+\gamma  p_i) dt +D dW_i,\\
 \dot {m_i} = (A + G)m_i +B\bar u +\gamma  p_i,\\
 \dot{ p_i} = -(A+G)^T  p_i   -(I-\Gamma)^T Q [\mathbb{E}x_i- (\Gamma m_i+\eta)],\\
dy_i=\left\{-A^T y_i +Q[x_i-(\Gamma m_i+\eta)]\right\} dt +\zeta_i dW_i,
 \end{cases}
\end{equation}
where $x_i(0)$ is given, $m_i(0)=m_0$, $ p_i(T)=H{\mathbb E} x_i(T)$, and $y_i(T)=-Hx_i(T)$.

As in the proof of Lemma \ref{lemma:uy}, we evaluate the first variation of $\bar J_i(u_i)$ at  $(x_i, m_i,  p_i, \check u_i)$ and can show $\delta \bar J_i=0$.
Since $\bar J_i$ is convex, this zero first variation condition implies that $\check u_i$ is an optimal control of ({P2}). By the same reasoning, $u_i'$ is also an optimal control. By strict convexity, we have $\check u_i=u_i'$. Subsequently, we have $(x_i, m_i,  p_i)=(x_i', m_i',  p_i')$ by Theorem \ref{theorem:p1}.
This  further implies $(y_i,\zeta_i)=(y_i', \zeta_i')$.~\qed

\section{The Solution of the Robust Game}
\label{sec:cons}

Note that Theorem \ref{theorem:ror} determines the strategy of a representative agent when $\bar u$ is fixed. Denote
\begin{align}\label{Nxm}
  x^{(N)} = \frac{1}{N}\sum_{i=1}^N x_i,
\quad y^{(N)} = \frac{1}{N}\sum_{i=1}^N y_i,\quad
 m^{(N)}=\frac{1}{N}\sum_{i=1}^N m_i ,\quad p^{(N)}=\frac{1}{N}\sum_{i=1}^N p_i.
 \end{align}
By \eqref{fbsdeubar}, we obtain
\begin{equation} \label{mfgeN}
\begin{cases}
dx^{(N)} =\left( Ax^{(N)} +BR^{-1}B^Ty^{(N)} +G m^{(N)}+\gamma  p^{(N)}\right) dt + \frac{D}{N}\sum_{i=1}^N dW_i,\\
 \frac{d m^{(N)}}{dt} = (A + G) m^{(N)} +B\bar u +\gamma  p^{(N)},\\
 \frac{d p^{(N)}}{dt} = -(A+G)^T p^{(N)}   -(I-\Gamma)^T Q \left[\mathbb{E}x^{(N)}- (\Gamma  m^{(N)}+\eta)\right] ,\\
dy^{(N)}=\left\{-A^T y^{(N)} +Q[x^{(N)}-(\Gamma m^{(N)}+\eta)]\right\} dt +\frac{1}{N}\sum_{i=1}^N\zeta_i dW_i,
\end{cases}
\end{equation}
where ${x}^{(N)}(0)= ({1}/{N}) \sum_{i=1}^N x_i(0)$, $m^{(N)}(0)=m_0$, $p^{(N)}(T)=H{\mathbb E} x^{(N)}(T)$, and $y^{(N)}(T)=-Hx^{(N)}(T)$.

As an approximation to \eqref{mfgeN}, we construct the following limiting system
\begin{equation} \label{mfge4}
\begin{cases}
\dot{ {\bf x}}= A {\bf x} +BR^{-1} B^T {\bf y} +G{\bf m} +\gamma {\bf p},    \\
\dot{ {\bf m }}= (A + G){\bf  m} +B\bar u +\gamma {\bf  p},\\
 \dot{{\bf  p}} = -(A+G)^T {\bf p }   -(I-\Gamma)^T Q [{\bf x}- (\Gamma  {\bf m}+\eta)] ,\\
\dot{{\bf y}}=-A^T {\bf y} +Q[{\bf x}-(\Gamma {\bf m}+\eta)] ,
\end{cases}
\end{equation}
where ${\bf x}(0)=  {\bf m}(0)=m_0$, ${\bf p}(T)=H {\bf x}(T)$, and ${\bf y}(T)=-H{\bf x}(T)$. This is a two point boundary value problem.

Note that ${\bf y}$ is intended as an approximation of $y^{(N)}$ when $N\to \infty$.
The consistency requirement imposes
\begin{equation}
\bar u= R^{-1} B^T{\bf y}. \label{cc}
\end{equation}

Under the  condition \eqref{cc}, the first two equations in \eqref{mfge4} coincide to give ${\bf x}={\bf m}$ for all $t\in [0,T]$. Consequently, we eliminate the equation of ${\bf x}$ and  introduce the new  system
\begin{equation} \label{mfge}
\begin{cases}
 \dot {\bf m }= (A + G){\bf  m} +BR^{-1}B^T {\bf y} +\gamma {\bf  p},\\
 \dot{\bf  p} = -(A+G)^T {\bf p }  -(I-\Gamma)^T Q [{\bf m}- (\Gamma  {\bf m}+\eta)],\\
\dot{\bf y}=-A^T {\bf y} +Q[{\bf m}-(\Gamma {\bf m}+\eta)],
\end{cases}
\end{equation}
 where $  {\bf m}(0)=m_0$, ${\bf p}(T)=H {\bf m}(T)$, and ${\bf y}(T)=-H{\bf m}(T)$. This is still a two point boundary value problem.  The next corollary follows from Theorem \ref{theorem:ror}.

\begin{corollary} \label{theorem:mpyp2}  Assume \emph{(H1)}-\emph{(H2)}.
Suppose that \eqref{mfge} has a unique solution $({\bf m},{\bf p},{\bf y})\in C^1([0,T]; \mathbb{R}^{3n})$   and take $\bar u= R^{-1} B^T{\bf y}$ in \eqref{fbsdeubar}. Then
\eqref{fbsdeubar} has  a unique solution $(x_i,  m_i,  p_i, y_i, \zeta_i)\in S[0,T]$. \qed
\end{corollary}


\subsection{The special case of same initial conditions}

Consider the special case where all agents have the same initial condition $x_i(0)=m_0$ for all $i\ge 1$.
 The  FBSDE  \eqref{fbsdeubar} defines a mapping
$$
\Lambda( \bar u)= R^{-1} B^T {\mathbb E}y_i,
$$
where we take $\bar u\in C([0,T];\mathbb{R}^{n_1} )$. Clearly $R^{-1} B^T{\mathbb E}y_i$ is a continuous $\mathbb{R}^{n_1}$-valued function of $t\in [0,T]$.

By the consistency requirement $\bar u =\Lambda(\bar u)$, we set $\bar u=  R^{-1} B^T {\mathbb E}y_i $ in the second equation of \eqref{fbsdeubar}
to obtain the equation system of the mean field game:
\begin{equation} \label{mfgeu}
\begin{cases}
dx_i =( Ax_i +BR^{-1}B^Ty_i +G m_i+\gamma  p_i) dt +D dW_i,\\
 \dot{ m_i }= (A + G) m_i +BR^{-1} B^T {\mathbb E} y_i +\gamma  p_i,\\
 \dot{ p_i} = -(A+G)^T p_i  -(I-\Gamma)^T Q [\mathbb{E}x_i- (\Gamma  m_i+\eta)] ,\\
dy_i=\left\{-A^T y_i +Q[x_i-(\Gamma m_i+\eta)]\right\} dt +\zeta_i dW_i,
\end{cases}
\end{equation}
where ${ x}_i(0)=m_i(0)=m_0$, $p_i(T)=H{\mathbb E} x_i(T)$, and $y_i(T)=-Hx_i(T)$.

An interesting fact is that the existence and uniqueness of a solution to  \eqref{mfgeu} is completely determined by the ODE system \eqref{mfge}
  without further using (H1)-(H2).

\begin{theorem}
\eqref{mfgeu} has a unique solution $(x_i,  m_i,  p_i, y_i, \zeta_i)\in S[0,T]$ if and only if \eqref{mfge} has a unique solution.
\end{theorem}

\proof
By Lemma \ref{lemma:sderelate}, \eqref{mfgeu} has a unique solution if and only if the FBSDE \eqref{mfge2} has a unique solution.
By Lemma \ref{lemma:sdeode} and Lemma \ref{lemma:3to4}-(iii), the FBSDE    \eqref{mfge2} has a unique solution if and only if \eqref{mfge} has a unique solution.
The theorem follows. \qed

\subsection{Existence of a solution to \eqref{mfge}}

To study the existence and uniqueness of a solution to \eqref{mfge}, we use a fixed point approach and introduce the equation system
\begin{equation} \label{mfgode3h}
\begin{cases}
 \dot{\bf m} = (A + G){\bf m} +h +\gamma  {\bf p},\\
 \dot{\bf  p} = -(A+G)^T {\bf  p}   -\widehat Q {\bf m}  + (I-\Gamma)^T  Q\eta,\\
\dot {\bf y}=-A^T {\bf  y} +Q[{\bf m}-(\Gamma {\bf  m}+\eta)] ,
\end{cases}
\end{equation}
where $h\in C([0,T]; \mathbb{R}^n) $,   $m(0)=m_0$, $ {\bf p}(T)=H{\bf m}(T)$, and  $ {\bf y}(T)=-H {\bf m} (T)$.
The next lemma  identifies a sufficient condition for \eqref{mfgode3h} to have a unique solution for any $h\in C([0,T]; \mathbb{R}^n)$.

\begin{lemma}
Suppose that the Riccati equation
\begin{align}\label{ricK}
\dot{K}+K(A+G)+(A+G)^TK -\gamma K^2 -\widehat Q=0,\qquad
K(T)=-H
\end{align}
 has a unique solution on $[0,T]$.
 Then  \eqref{mfgode3h} defines a mapping
  from $C([0,T]; {\mathbb R}^n)$ to itself:
$$
\Lambda_1: h \longmapsto BR^{-1}B^T {\bf  y}.
$$

\end{lemma}

{\it Proof}.
We write ${\bf  p}=-K{\bf  m} +\phi$ for \eqref{mfgode3h} and obtain the ODE
$$
\dot \phi = -(A+G-\gamma K) \phi + Kh+ (I-\Gamma)^T Q\eta , \quad \phi(T)=0.
$$
It follows that
$$
\dot{\bf m}= (A+G-\gamma K)  {\bf m} +h +\gamma \phi.
$$
Let the fundamental solution matrices of the two ODEs
$$
\dot \varphi =(A+G-\gamma K)\varphi, \quad  \dot \psi =-(A+G-\gamma K)^T\psi
$$
 be $\Phi(t,s)$ and $\Psi(t,s)$, respectively,  with $\Phi(s,s)=\Psi(s,s)=I$.
Then $\Psi(t,s)= \Phi^T(s,t)$.
We obtain
$$
\phi(t)=-\int_t^T \Psi(t,s_1) [K(s_1)h(s_1)+ (I-\Gamma)^T Q\eta] ds_1.
$$
This in turn gives
\begin{align*}
{\bf m}(t)=\ & \Phi(t,0) m_0+ \int_0^t \Phi(t, s_1) h(s_1) ds_1\\
& -\gamma \int_0^t \Phi(t, s_2) \int_{s_2}^T \Psi(s_2,s_1)[K(s_1)h(s_1)+ (I-\Gamma)^T Q\eta] ds_1  ds_2 .
\end{align*}
We further solve
$$
{\bf y}(t)=  -\int_t^T e^{-A^T(t-s_3)}Q[(I-\Gamma){\bf  m} (s_3)-\eta] ds_3 - e^{-A^T(t-T)} H {\bf m}(T),
$$
which implies ${\bf y}\in C([0,T]; {\mathbb R}^n)$. The lemma follows. \qed

To simplify the existence analysis for \eqref{mfge} in this section, we consider the case $H=0$.
Below  $\Upsilon_k$  denotes a continuous function of $t$ which does not depend on $h$ and  can be easily determined.
Consequently,
\begin{align*}
{\bf y}(t)&= -\int_t^T e^{-A^T(t-s_3)}Q[(I-\Gamma){\bf  m} (s_3)-\eta] ds_3\\
&= -\int_t^T e^{-A^T(t-s_3)}Q(I-\Gamma){\bf  m} (s_3) ds_3 +\Upsilon_1(t)\\
&=  -\int_t^T e^{-A^T(t-s_2)}Q(I-\Gamma)  \int_0^{s_2} \Phi(s_{2}, s_1) h(s_1) ds_1 ds_2\\
&\quad +  \gamma \int_t^T e^{-A^T(t-s_3)}Q(I-\Gamma)   \int_0^{s_3} \Phi(s_3, s_2) \int_{s_2}^T \Psi(s_2,s_1)K(s_1)h(s_1) ds_1  ds_2  ds_3\\
&\quad +\Upsilon_2(t).
\end{align*}
Now we have
\begin{align*}
\Lambda_1(h)(t)=  &\ BR^{-1}B^T {\bf  y }(t) \\
=& -BR^{-1} B^T  \int_t^T e^{-A^T(t-s_2)}Q(I-\Gamma)  \int_0^{s_2} \Phi(s_{2}, s_1) h(s_1) ds_1 ds_2\\
& + \gamma BR^{-1} B^T \int_t^T e^{-A^T(t-s_3)}Q(I-\Gamma)   \int_0^{s_3} \Phi(s_3, s_2) \int_{s_2}^T \Psi(s_2,s_1)K(s_1)h(s_1) ds_1  ds_2  ds_3\\
& +BR^{-1} B^T\Upsilon_2(t)\\
\eqdef &\ \Lambda_0 (h)(t)  +BR^{-1} B^T\Upsilon_2(t).
\end{align*}
It is clear that $\Lambda_0$ is from $C([0,T]; {\mathbb R}^n)$ to itself. \qed

Define the constants
\begin{align*}
&  c_1= \max_{t\in [0,T]} |K(t)|,\quad c_2= \max_{0\le t,s \le T}|\Phi (t,s)|, \\
&c_3= \max_{t\in [0,T]} \int_t^T|e^{A (s-t)}|sds, \quad   c_4= \max_{t\in [0,T]}\int_t^T |e^{A(s-t)}| (Ts-\frac{s^2}{2}) ds.
\end{align*}
 Note that $Ts-\frac{s^2}{2}\ge 0$ for $s\in [0,T]$.
 Denote $|h|=\max_{t\in[0,T]} |h(t)|$.

\begin{theorem} \label{theorem:c1c4}
Assume $H=0$.  If
\begin{equation}
c_2|BR^{-1}B^T|\cdot |Q(I-\Gamma)|  \left(c_3+\gamma c_{1}c_2c_4\right)<1, \label{less1}
\end{equation}
 then   \eqref{mfge} has a unique solution.
\end{theorem}

{\it Proof.}
 For each $t$,
\begin{align*}
|\Lambda_0 (h)(t) | \le & \ c_2|h|\cdot |BR^{-1} B^T|\cdot |Q(I-\Gamma)|
\int_t^T |e^{A^T(s_2-t)}| s_2 ds_2  \\
&+ \gamma c_1c_2^2|h|\cdot|BR^{-1} B^T|\cdot |Q(I-\Gamma)|
\int_t^T |e^{A^T(s_3-t)}|  \int_0^{s_3} \int_{s_2}^{T} ds_1 ds_2 ds_3\\
=&\ c_2|h|\cdot |BR^{-1} B^T|\cdot |Q(I-\Gamma)|
\int_t^T |e^{A(s_2-t)}| s_2 ds_2\\
&+\gamma c_1c_2^2|h|\cdot|BR^{-1} B^T|\cdot |Q(I-\Gamma)|
\int_t^T |e^{A(s_{3}-t)}|  (Ts_3-\frac{s_3^2}{2})  ds_{3}\\
=&\ c_2|BR^{-1}B^T|\cdot |Q(I-\Gamma)|  \left(c_3+\gamma c_{1}c_2
c_4\right).
\end{align*}
Hence, $\Lambda_1$ is a contraction and   has a  unique fixed point.
So \eqref{mfge} has a unique solution.
 \qed

The  constants $c_1, \ldots, c_4$ in \eqref{less1} do not depend on $BR^{-1}B^T$. If $BR^{-1}B^T$ is suitably small, \eqref{less1} can be ensured.

\begin{example}
Consider the system with parameters given by Example \ref{example:A}. Take $T=1.3$. In analogue  to \eqref{numP}, we can solve $K(t)$ on $[0,T]$ for \eqref{ricK}. It can be shown that $K(t)\le 0$ for $t\in [0,T]$ and $|K(t)|$ attains its maximum on $[0,T]$ at $t=0$. We have $K(0)= -0.171417$
which gives $c_1=0.171417$. So $c_2\le e^{(A+G+|K(0)|)T}=3.312961$.
Furthermore,
$$
c_3\le \int_0^Te^{As}s ds= 1.318243 , \quad c_4\le
\int_0^Te^{As} (Ts-\frac{s^2}{2})  ds =1.112937.
$$
Subsequently,
$$
c_2|BR^{-1}B^T|\cdot |Q(I-\Gamma)|  \left(c_3+\gamma c_{1}c_2c_4\right)
\le 0.861493.
$$
So \eqref{less1} holds.
\end{example}

%


\begin{remark}
\emph{For the two-point boundary value problem, the
contraction estimate in
 the fixed point method may be  conservative  and typically works on small time intervals for  the solvability of \eqref{mfge}  (see, e.g., Ch.1, Sec. 5, \cite{MY99}). }
\end{remark}

 We continue to derive another condition under which \eqref{mfge} is solvable without restriction to a small time horizon. To this end, we first rewrite \eqref{mfge} in the following form:
\begin{equation}\label{E01}
\left\{ \begin{array}{lll}
\left(
       \begin{array}{c}
        \dot{m} \\
        \dot{p} \\
         \dot{y} \\
       \end{array}
     \right)
=
\widetilde{A}\left(
       \begin{array}{c}
         m \\
         p \\
         y \\
       \end{array}
     \right)+\widetilde\eta,\\
m(0)=m_0,\quad p(T)=Hm(T),\quad y(T)=-Hm(T),
\end{array}\right.
\end{equation}
where
$$
\widetilde{A}=\left(
            \begin{array}{ccc}
              A+G & \gamma & BR^{-1}B^{T}\\
             -(I-\Gamma)^{T}Q(I-\Gamma) & -(A+G)^{T} & 0 \\
              Q(1-\Gamma) & 0 & -A^{T} \\
            \end{array}
          \right),
\quad
\widetilde \eta=\left(
       \begin{array}{c}
         0 \\
        (I-\Gamma)^{T}Q\eta \\
         -Q\eta \\
       \end{array}
     \right).
$$Then, by the variation of constant formula, we have
\begin{equation}\label{E02}
\left(
       \begin{array}{c}
         m(t) \\
         p(t) \\
         y(t) \\
       \end{array}
\right)=\Theta(t)\left(
       \begin{array}{c}
         m_0 \\
        \mu \\
        \nu \\
       \end{array}
\right)+\Theta(t)\int_0^t \Theta^{-1}(s)\widetilde \eta ds,
\end{equation}
where $\Theta(t)=e^{\widetilde{A}t}$ and $p,y$ have the initial conditions $p(0)=\mu$, $y(0)=\nu$. Noting the terminal condition in \eqref{E01}, now we present the following result.
\begin{proposition} \emph{\cite[Ch.2, Sec. 3]{MY99}}
If for given $T>0,$ $\det(\widetilde{\Theta}(T)) \neq 0,$ where$$
\widetilde{\Theta}(T)=\left(
                           \begin{array}{ccc}
                             -H & I & 0 \\
                             H & 0 & I \\
                           \end{array}
                         \right)\Theta(T)\left(
                                           \begin{array}{cc}
                                             0 & 0 \\
                                             I & 0 \\
                                             0 & I \\
                                           \end{array}
                                         \right),
$$then \eqref{mfge} has a unique solution on $[0, T]$ for any initial value $m_0.$
\end{proposition}


 For illustration, we give the following example.

\begin{example} \label{example:A3}
Consider system \eqref{x1}-\eqref{jcost} with all parameters being scalar-valued and $\Gamma=1$, $H=0$.
 We calculate
 $$  \mathcal{A}=\left(\begin{array}{cc} A+G & -\gamma  \\
0 & -(A+G)\\ \end{array}\right), \quad
 \widetilde{A}=\left(
            \begin{array}{ccc}
             A+G & \gamma & R^{-1}B^{2}\\
             0 & -(A+G) & 0 \\
             0 & 0 & -A \\
            \end{array}
          \right),$$
where ${\mathcal A}$ is defined in
Proposition \ref{lemma:ricj}.
By direct computations, we obtain
$$
\det\{[(0, I)e^{\mathcal{A}t}(0, I)^{T}]\}=e^{-(A+G)t}>0
$$
for all $t\in [0,T]$, which ensures \emph{(H1)} by
Proposition \ref{lemma:ricj}.
 Moreover, $b_3=0$ gives $C_q=0$ in \eqref{Cqnu} so that \emph{(H2)} always holds true for $R>0$.                
Finally, $\det(\widetilde{\Theta}(t))=e^{-(2A+G)t}>0,$ and
subsequently,  \eqref{mfge} has a unique solution on any interval
$[0, T].$ To summarize, \emph{(H1)}, \emph{(H2)} and the solvability of \eqref{mfge}
are all satisfied by the system.
\end{example}

Note that the solvability of \eqref{mfge} in Example  \ref{example:A3} does not depend on the value of $R^{-1}B^{2}$, which is different from the condition in Theorem \ref{theorem:c1c4}.

\section{Error Estimate of the Mean Field Approximation}
\label{sec:error}

 We suppose that \eqref{mfge} has  a unique solution
$({\bf m},{\bf p},{\bf y})$ and accordingly take $\bar u$ in \eqref{fbsdeubar} as
\begin{align}
\bar u^* = R^{-1} B^T {\bf y}. \label{uby2}
\end{align}
The FBSDE system \eqref{fbsdeubar} now becomes
\begin{equation}\label{fbsdeubar2}
\begin{cases}
dx_i =( Ax_i +BR^{-1}B^Ty_i +G m_i+\gamma  p_i) dt +D dW_i,\\
 \dot{m_i} = (A + G) m_i +B\bar u^* +\gamma  p_i,\\
 \dot{p_i} = -(A+G)^T  p_i   -(I-\Gamma)^T Q [\mathbb{E}x_i- (\Gamma m_i+\eta)] ,\\
dy_i=\left\{-A^T y_i +Q[x_i-(\Gamma m_i+\eta)]\right\} dt +\zeta_i dW_i,
 \end{cases}
\end{equation}
where $x_i(0)$ is given, $m_i(0)=m_0$, $ p_i(T)=H{\mathbb E} x_i(T)$, and $y_i(T)=-Hx_i(T)$.
By Corollary \ref{theorem:mpyp2}, this FBSDE has a unique solution.
In the game of $N$ players,
let $y_i$ be solved from \eqref{fbsdeubar2} and  denote the control for ${\cal A}_i$ by
\begin{align}
\hat u_i= R^{-1}B y_i, \qquad 1\le i\le N,  \label{huiN}
\end{align}
which is a well defined process in $L_{\cal F}^2(0,T;\mathbb{R}^{n_1})$.

For $\hat u^{(N)}=(1/N)\sum_{i=1}^N \hat u_i $,
we aim to estimate
$$
{\mathbb E}|\hat u^{(N)}(t)-\bar u^*(t)|^2.
$$
Note that $\hat u_1, \ldots, \hat u_N$ are independent, but they are not necessarily with the same distribution due to possibly different initial states of the agents. This fact will somehow complicate our error estimate. The key result of this section is the following theorem.

\begin{theorem}
 \label{theorem:hubu}
Assume that \emph{(H1)}-\emph{(H2)} hold and that \eqref{mfge} has a unique solution.  We have
$$\sup_{0\le t\le T}{\mathbb E}|\hat{u}^{(N)} -\bar u^*|^2=O(1/N)+O(|x^{(N)}(0)-m_0|^2), $$
where $x^{(N)}(0)= (1/N) \sum_{i=1}^N x_i(0)$. \qed
\end{theorem}

The proof of Theorem \ref{theorem:hubu} is provided in the remaining part of this section. To do this, we need to prove some lemmas under the assumption of the theorem. Recalling \eqref{Nxm}, we  take $\bar u=\bar u^*$ in \eqref{mfgeN}  to write
\begin{equation} \label{mfgeN01}
\begin{cases}
dx^{(N)} =\left( Ax^{(N)} +BR^{-1}B^Ty^{(N)} +G m^{(N)}+\gamma  p^{(N)}\right) dt + \frac{D}{N}\sum_{i=1}^N dW_i,\\
 \frac{d m^{(N)}}{dt} = (A + G) m^{(N)} +B\bar u^* +\gamma  p^{(N)},\\
 \frac{d p^{(N)}}{dt} = -(A+G)^T p^{(N)}   -(I-\Gamma)^T Q \left[\mathbb{E}x^{(N)}- (\Gamma  m^{(N)}+\eta)\right] ,\\
dy^{(N)}=\left\{-A^T y^{(N)} +Q[x^{(N)}-(\Gamma m^{(N)}+\eta)]\right\} dt +\frac{1}{N}\sum_{i=1}^N\zeta_i dW_i,
\end{cases}
\end{equation}
where ${x}^{(N)}(0)= ({1}/{N}) \sum_{i=1}^N x_i(0)$, $m^{(N)}(0)=m_0$, $p^{(N)}(T)=H{\mathbb E} x^{(N)}(T)$, and $y^{(N)}(T)=-Hx^{(N)}(T)$.

Denote the ODE system
\begin{equation} \label{mfge4N}
\begin{cases}
\dot{\bf x}_N= A {\bf x}_N +BR^{-1} B^T {\bf y}_N +G{\bf m}_N +\gamma {\bf p}_N,    \\
 \dot{\bf m }_N= (A + G){\bf  m}_N +B\bar u^* +\gamma {\bf  p}_N,\\
 \dot{\bf  p}_N = -(A+G)^T {\bf p }_N   -(I-\Gamma)^T Q [{\bf x}_N- (\Gamma  {\bf m}_N+\eta)] ,\\
\dot{\bf y}_N=-A^T {\bf y}_N +Q[{\bf x}_N-(\Gamma {\bf m}_N+\eta)],
\end{cases}
\end{equation}
where ${\bf x}_N(0) = ({1}/{N})\sum_{i=1}^N x_i(0)$,  $  {\bf m}_N(0)=m_0$, ${\bf p}_N(T)=H {\bf x}_N(T)$, and ${\bf y}_N(T)=-H{\bf x}_N(T)$. The initial condition ${\bf x}_N(0)$ is  different from that of \eqref{mfge4}.

\begin{lemma} \label{lemma:eN}
\eqref{mfge4N} has a unique solution which can be denoted as
$$
({\bf x}_N, {\bf m}_N, {\bf p}_N, {\bf y}_N)=({\mathbb E}x^{(N)}, m^{(N)}, p^{(N)}, {\mathbb E}y^{(N)}).
$$
  \end{lemma}

 \proof Existence follows by taking expectation in \eqref{mfgeN01}.
 To show uniqueness, suppose that \eqref{mfge4N} has two different solutions $({\bf x}_N, {\bf m}_N, {\bf p}_N, {\bf y}_N)$ and $({\bf x}_N', {\bf m}_N', {\bf p}_N', {\bf y}_N')$. Then for any $\lambda\in \mathbb{R}$,
$$
(x_i,  m_i, p_i, y_i, \zeta_i) +\lambda ({\bf x}_N-{\bf x}_N', {\bf m}_N-{\bf m}_N', {\bf p}_N-{\bf p}_N', {\bf y}_N-{\bf y}_N',0)
$$
satisfies \eqref{fbsdeubar}, which is a contradiction to Theorem
\ref{theorem:ror}. Uniqueness follows.  \qed

\begin{lemma}\label{lemma:xyN}
We have
$$
\sup_{0\le t\le T}\left( {\mathbb E}|x^{(N)} -{\mathbb E} x^{(N)}|^2 +{\mathbb E}|y^{(N)}-{\mathbb E} y^{(N)} |^2\right) =O(1/N).  $$

\end{lemma}

\proof
 Define
$$(\theta_1, \theta_2)= (x^{(N)} -{\mathbb E} x^{(N)} ,y^{(N)}-{\mathbb E} y^{(N)} ). $$
By \eqref{mfgeN01}, \eqref{mfge4N} and Lemma \ref{lemma:eN},
\begin{equation*}
\begin{cases}
d\theta_1= (A \theta_1+BR^{-1}B^T \theta_2)dt +\frac{D}{N}\sum_{i=1}^N dW_i,\\
d\theta_2= (-A^T \theta_2+Q \theta_1)dt  +\frac{1}{N }\sum_{i=1}^N \zeta_i dW_i,
\end{cases}
\end{equation*}
where $\theta_1(0)=0$ and $\theta_2(T)= -H \theta_1(T)$.

Let $P$ be the solution of the Riccati equation
\begin{equation*}
\dot P+A^TP+PA-PBR^{-1}B^TP+Q=0, \quad P(T)=H. 
\end{equation*}
 Denote
$
\theta_2 = -P\theta_1+\psi,
$
where $\psi(T)=0$.
This gives the equation
\begin{align}
d\psi = -(A-BR^{-1} B^T P)^T \psi dt +\frac{1}{N}\sum_{i=1}^N (PD+\zeta_i)dW_i,\nonumber
\end{align}
where $\psi(T)= 0$. There is a unique solution $\psi=0$ for $t\in [0,T]$. This implies
\begin{align}
d\theta_1= (A -BR^{-1}B^T P)\theta_1dt +\frac{D}{N}\sum_{i=1}^N dW_i.\nonumber
\end{align}
Hence,
$
\sup_{0\le t\le T}{\mathbb E}|\theta_1(t)|^2= O(1/N).
$
The lemma follows since $\theta_2=-P\theta_1$. \qed

When $({\bf m},{\bf p},{\bf y})$ is a unique solution of \eqref{mfge}, it can be shown that  $({\bf x},{\bf m},{\bf y},{\bf p}) \defeq ({\bf m},{\bf m},{\bf y},{\bf p})$ is a unique solution of \eqref{mfge4} under the condition \eqref{cc}.


\begin{lemma} \label{lemma:xmpy0}
We have
\begin{align*}
\sup_{0\le t\le T} [|{\bf x}_N- {\bf x} | +|{\bf m}_N- {\bf m} |+ |{\bf p}_N- {\bf p} | + |{\bf y}_N- {\bf y} |]= O(|x^{(N)}(0)-m_0|).
\end{align*}
\end{lemma}

\proof
Consider
\begin{equation} \label{h14}
\begin{cases}
\dot { h}_1= A { h_1} +BR^{-1} B^T { h_4} +G{ h_2} +\gamma { h_3},    \\
 \dot { h}_2 = (A + G){h_2} +\gamma {h_3},\\
 \dot{h}_3 = -(A+G)^T {h_3 }   -(I-\Gamma)^T Q ({h_1}- \Gamma  {h_2}) ,\\
\dot{h}_4=-A^T {h_4} +Q({h_1}-\Gamma {h_2}),
\end{cases}
\end{equation}
where ${ h_1}(0)$ is given,  ${h_2}(0)=0$, ${ h_3}(T)=H {h_1}(T)$, and ${ h_4}(T)=-H{h_1}(T)$. It is constructed as a homogeneous version of \eqref{mfge4N}.
 We claim that \eqref{h14} has a unique solution for any given value of $h_1(0)$.
If this were not true, there would exist $h(0)$
such that \eqref{h14} has multiple solutions which, in turn, can be used to construct multiple solutions  to  \eqref{mfge4N}. This would give a contradiction to Lemma \ref{lemma:eN}.

It is clear that
$$
({\bf x}_N-{\bf x},{\bf m}_N-{\bf m}, {\bf p}_N-{\bf p}, {\bf y}_N-{\bf y}  )\eqdef(h_1, h_2, h_3, h_4),
$$
is a solution of \eqref{h14} with $h_1(0)= {\bf x}_N(0)-m_0$.

Let $e_1, \ldots e_n$ be a canonical basis of $\mathbb{R}^n$. For $h_1(0)=e_k$, we obtain a solution of \eqref{h14}, denoted by $h^k=(h_1^k, h_2^k, h_3^k, h_4^k)$. Let $(z)_k$ be the $k$th component of a vector $z$. We may uniquely denote $ ({\bf x}_N-{\bf x},{\bf m}_N-{\bf m}, {\bf p}_N-{\bf p}, {\bf y}_N-{\bf y}  )$ as a linear combination of $h^1, \ldots, h^n$:
$$
({\bf x}_N-{\bf x},{\bf m}_N-{\bf m}, {\bf p}_N-{\bf p}, {\bf y}_N-{\bf y}  )=\sum_{k=1}^n (x_N(0)-m_0)_k (h_1^k, h_2^k, h_3^k, h_4^k).
$$
The lemma follows readily.\qed

{\bf Proof of Theorem \ref{theorem:hubu}.}
For $\bar u=\bar u^*$, we write $\hat u^{(N)}= R^{-1} B^T (1/N)\sum_{i=1}^Ny_i= R^{-1}B^T y^{(N)}$.
We have
\begin{align*}
|\hat u^{(N)}-\bar u^* |^2 &= {\mathbb E}|R^{-1}B^T(y^{(N)}- {\bf y})|^2 \\
   &\le C {\mathbb E}|y^{(N)}-{\bf y}|^2 \\
   &= CE|y^{(N)} -{\mathbb E} y^{(N)}+{\mathbb E}y^{(N)}-{\bf y}|^2\\
   &\le C (1/N) +C |{\bf y}_N-{\bf y}|^2\\
   &=O(1/N)+ O(|x^{(N)}(0)-m_0|^2).
\end{align*}
The second inequality follows from Lemmas \ref{lemma:eN} and \ref{lemma:xyN}, and the last step follows from Lemma \ref{lemma:xmpy0}.
 \qed

\section{Robust Nash Equilibrium}
\label{sec:nash}

Throughout this section,   we assume that \eqref{mfge} has a unique solution and take $\bar u=\bar u^*$ determined by
\eqref{uby2}.
For $f\in L^2(0,T; \mathbb{R}^n)$ and  $u_i\in L_{\cal F}^2(0,T; \mathbb{R}^{n_1})$, $1\le i \le N$, recall the worst case cost
\begin{align*}
J^{\rm wo}_i(u_i, u_{-i})= \sup_{f\in L^2(0,T;\mathbb{R}^n) }J_i(u_i, u_{-i},f).
\end{align*}
It is clear that for each $i$ and any $(u_i, u_{-i})$, $\sup_f J_i( u_i,  u_{-i},f)\ge 0$.

Consider the set of strategies $(\hat u_i, \hat u_{-i})$  given by \eqref{huiN} for a population of $N$ players with dynamics \eqref{x1}.
 It should be emphasized that we only use \eqref{fbsdeubar2}-\eqref{huiN} to make a well defined process $\hat u_i$ in $L_{\cal F}^2(0,T; \mathbb{R}^{n_1})$ which should not be understood as a feedback strategy.
 The main result of this section is the next theorem which characterizes the performance of this set of strategies.

\begin{theorem}
\label{theorem:nash}  Assume
\emph{(i)} \emph{(H1)}-\emph{(H2)} hold;
\emph{(ii)} $\sup_{i\ge 0} |x_i(0)|\le M_0$ where $M_0$ does not depend on $N$; \emph{(iii)} \eqref{mfge} has a unique solution. Then the set of strategies $(\hat u_1, \ldots, \hat u_N)$ given by \eqref{huiN} is a robust $\varepsilon_N$-Nash equilibrium for the $N$ players, i.e.,
\begin{align}
J_i^{\rm wo}(\hat u_i, \hat u_{-i})-\varepsilon_N \le  \inf_{u_i\in {\cal U}} J_i^{\rm wo}(u_i, \hat u_{-i})  \le J_i^{\rm wo}(\hat u_i, \hat u_{-i}),
\label{jwo}
\end{align}
where $0\le \varepsilon_N =O(1/\sqrt{N}+ |x^{(N)}(0)-m_0|)$ and $x^{(N)}(0)=(1/N)\sum_{j=1}^N x_j(0)$. \qed
\end{theorem}

The rest part of this section is devoted to the proof of Theorem \ref{theorem:nash}.
 For any given
 $f\in L^2(0,T;\mathbb{R}^n)$,
 denote the  state processes of \eqref{x1} corresponding to $(\hat u_i, \hat u_{-i}, f)$    by $\hat x_j$, $1\le j \le N$, and $\hat x^{(N)}= (1/N)\sum_{j=1}^N \hat x_j$.
Denote
\begin{equation}
\dot{\bar m} =(A+G) \bar m +B\bar u^*+f ,  \quad  \bar m(0)=m_0  \label{xbf}
\end{equation}

 All subsequent  lemmas are proved under the assumptions of
 Theorem \ref{theorem:nash}.

\begin{lemma} \label{lemma:hx-bx}
We have
\begin{align}
\sup_{0\le t\le T, f} {\mathbb E}|\hat x^{(N)}
-\bar m|^2 \le C(1/N +| x^{(N)}(0)-m_0|^2). \nonumber
\end{align}
\end{lemma}

\proof Note that
$$
d \hat x^{(N)} = [(A+G) \hat x^{(N)} +B\hat u^{(N)} +f ]dt +(D/N)\sum_{i=1}^N dW_i.
$$
 Therefore,
\begin{align}
d (\hat x^{(N)}-\bar m) = [(A+G) (\hat x^{(N)}-\bar m) +B(\hat u^{(N)}-\bar u^*) ]dt +(D/N)\sum_{i=1}^N dW_i. \nonumber
\end{align}
By linear SDE estimates,
\begin{align*}
{\mathbb E}|\hat x^{(N)}(t)-\bar m(t) |^2 \le\ &C |x^{(N)}(0)-m_0|^2 +C/N\\
&+  C{\mathbb E} \int_0^t  |\hat u^{(N)}(\tau)-\bar u^*(\tau)|^2d\tau.
\end{align*}
By Theorem \ref{theorem:hubu}, the lemma follows.
\qed

\begin{lemma}
There exists a constant $\hat C_0$ independent of $N$ such that
$$
\max_{1\le i\le N} \sup_f J_i(\hat u_i, \hat u_{-i},f)\le \hat C_0.
$$
\end{lemma}

{\it Proof.}
Denote
\begin{align} \label{xip}
dx_i' = (Ax'_i +B\hat u_i +G\bar m+f) dt +DdW_i,
\end{align}
where $x_i'(0)=x_i(0)$.
By Lemma \ref{lemma:hx-bx}, it is easy to show  $$\sup_{0\le t\le T,f}{\mathbb E}|\hat x_i(t)-x'_i(t)|^2\le C(1/N +| x^{(N)}(0)-m_0|^2).$$
We have
\begin{align}
J_i(\hat u_i, \hat u_{-i},f)\le\ & \bar J_i(\hat u_i, f) + {\mathbb  E}\int_0^T  |(\hat x_i-x'_i)+\Gamma(\bar m-\hat x^{(N)})|_Q^2dt + {\mathbb E} |\hat x_i(T)-x_i'(T)|_H^2 \nonumber \\
& +2 {\mathbb  E}\int_0^T[x'_i -(\Gamma \bar m+\eta)]^TQ [(\hat x_i-x'_i)+\Gamma(\bar m-\hat x^{(N)})]dt \nonumber \\
& + 2{\mathbb E}[x_i'^{T}(T) H (\hat x_i(T)-x_i'(T))]. \label{jNjb}
 \end{align}
Combining Lemma \ref{lemma:coer} with condition (ii)
in Theorem \ref{theorem:nash}, we obtain
  \begin{align}
 \bar J_i(\hat u_i,f)\le C -(\epsilon_0/2) \|f\|^2_{L^2} \label{je0}
\end{align}
for $\epsilon_0>0$, where $C$ does not depend on $(i,N)$.
Since neither $\hat x_i-x_i'$ nor $\bar m-\hat x^{(N)}$ depend on $f$,
there exists a constant $C_1$ such that
\begin{align}
&\left|{\mathbb  E}\int_0^T[x'_i -(\Gamma \bar m+\eta)]^TQ [(\hat x_i-x'_i)+\Gamma(\bar m-\hat x^{(N)})]dt\right|  \nonumber \\
&\le C_1\left( {\mathbb  E}\int_0^T|x'_i -(\Gamma \bar m+\eta)|^2_Qdt \right)^{1/2}
\nonumber\\
&\le C_2 (1+\|f\|^2_{L^2})^{1/2} \nonumber \\
&\le C_3 +(\epsilon_0/16)\|f\|^2_{L^2}, \label{xmf}
\end{align}
where the second inequality follows from elementary estimates based on the solutions of \eqref{xbf} and \eqref{xip}. Similarly,
\begin{align}
{\mathbb E}[x_i'^{T}(T) H (\hat x_i(T)-x_i'(T))]\le C_4
+(\epsilon_0/16)\|f\|^2_{L^2}. \label{xhxh}
\end{align}
 Finally combining \eqref{jNjb}-\eqref{xhxh} with Lemma \ref{lemma:hx-bx} leads to
 $$
 J_i(\hat u_i, \hat u_{-i}, f) \le C-(\epsilon_0/4) \|f\|^2_{L^2}.
 $$
 The lemma follows.
 \qed

Consider the set of strategies $(u_i, \hat u_{-i})$ and the corresponding state processes
\begin{align}
&dx_i = (Ax_i +Bu_i + Gx^{(N)} +f) dt +D dW_i, \label{xiuihu}\\
&dx_j=( Ax_j +B\hat u_j +Gx^{(N)} + f) dt +D dW_j, \quad 1\le j\le N,\ j\ne i. \label{xjhu}
\end{align}

\begin{lemma}
If $u_i$ in \eqref{xiuihu} satisfies $\sup_f J_i( u_i, \hat u_{-i},f)\le \hat C_0$,  there exists $\hat C_1$ independent of $N$ such that
\begin{align}  \label{uic1}
{\mathbb E}\int_0^T |u_i(t)|^2 dt \le \hat C_1.
\end{align}
\end{lemma}

{\it Proof.}
Suppose
$
\sup_fJ_i(u_i, \hat u_{-i},f)\le \hat C_0.
$
Then for any $f$,
$$
\mathbb{E}\int_0^T \left(|x_i- (\Gamma x^{(N)}+\eta)|_Q^2
 + u_i^T R u_i -\frac{1}{\gamma} |f(t)|^2\right) dt +{\mathbb E} [x_i^T(T) Hx_i(T)] \le \hat C_0,
$$
where $(x_1, \cdots, x_N)$ is generated by $(u_i, \hat u_{-i})$ and $f$.  Taking $f=0$, we obtain
$$
\mathbb{E}\int_0^T \left(|x_i- (\Gamma x^{(N)}+\eta)|_Q^2
 + u_i^T R u_i \right) dt\le \hat C_0.
$$
Therefore, \eqref{uic1} holds. \qed

Let ${\cal U}_{\hat C_1}$ denote the set of processes $u_i\in L_{\cal F}^2(0,T; \mathbb{R}^{n_1})$ which  satisfy \eqref{uic1}.
For \eqref{xiuihu}-\eqref{xjhu}, denote $x^{(N)}= (1/N)\sum_{j=1}^N x_j$.

\begin{lemma} \label{lemma:xxbar}
Suppose $u_i\in {\cal U}_{\hat C_1}$ in \eqref{xiuihu}. Then
$$
\sup_{0\le t\le T, f,u_i\in {\cal U}_{\hat C_1}}{\mathbb E}| x^{(N)}(t)-\bar m(t) |^2= O(1/N +|x^{(N)}(0)-m_0|^2  ).
$$
\end{lemma}

\proof
Rewrite  \eqref{xiuihu} in the form
\begin{align}
dx_i= [Ax_i+B\hat u_i +Gx^{(N)} +f] dt +B(u_i-\hat u_i) dt +DdW_i.  \label{xiunew}
\end{align}
By \eqref{xjhu} and \eqref{xiunew},
\begin{align}
dx^{(N)} = [(A+G)x^{(N)} +B \hat u^{(N)}+f ] dt +\frac{B}{N}(u_i-\hat u_i)dt +\frac{D}{N} \sum_{j=1}^N dW_j,\nonumber
\end{align}
which combined with \eqref{xbf} gives
\begin{align*}
d(x^{(N)}-\bar m)= \ &[ (A+G) (x^{(N)}-\bar m)+B(\hat u^{(N)}-\bar u^*)]dt \\
&+ \frac{B}{N}  (u_i-\hat u_i) dt +\frac{D}{N} \sum_{j=1}^N dW_j.
\end{align*}
By Theorem \ref{theorem:hubu} and the fact
 ${\mathbb E}\int_0^T|u_i-\hat u_i|^2\le C$ for all $u_i\in {\cal U}_{\hat C_1}$, where the constants $C$ do not depend on $(f,u_i)$, elementary SDE estimates lead  to
$$
\sup_{0\le t\le T,f} {\mathbb E}|x^{(N)}(t)-\bar m(t)|^2\le C(1/N+|x^{(N)}(0)-m_0|^2),
$$
where $C$ does not depend on $u_i$.
The lemma follows.
\qed

\begin{lemma}\label{lemma:jjb}
For each $u_i\in {\cal U}_{\hat C_1}$, $\sup_fJ_i(u_i, \hat u_{-i},f)  $ is finite and attained by some $f$ depending on $u_i$ and so denoted as $f_{u_i}$. Moreover,
$$
\sup_{u_i\in {\cal U}_{\hat C_1}}|\sup_fJ_i(u_i, \hat u_{-i},f) - \bar J_i (u_i,  \hat f_{u_i})|=O(1/\sqrt{N}+|x^{(N)}(0)-m(0)|),
$$
where $\hat f_{u_i}$ is determined by Theorem \ref{theorem:p1} for the given $u_i$.
\end{lemma}

\proof Note that we have
\begin{align}
&dx_i=[ Ax_i+Bu_i+ G\bar m + G(x^{(N)}-\bar m) +f ] dt+DdW_i, \label{xier}\\
& \dot{\bar m} = (A+G) \bar m + B\bar u^*+f , \nonumber
\end{align}
where $\bar m(0)=m_0$. Define the auxiliary process
\begin{align}
dx_i^\dag = (Ax_i^\dag+Bu_i+ G\bar m +  f ) dt+DdW_i,\nonumber
\end{align}
where $x_i^\dag(0)=x_i(0)$ and  $(u_i, f, W_i)$ is the same as in \eqref{xier}.
By Lemma \ref{lemma:xxbar}, it is easy to show
\begin{align}
\sup_{0\le t\le T, f} {\mathbb E} |x_i(t)-x_i^\dag(t)|^2 =O(1/N+ |x^{(N)}(0)-m(0)|^2  ).\label{xxdag}
\end{align}

We have the relation
\begin{align*}
|x_i -(\Gamma x^{(N)} +\eta) |_Q^2
=\ & |x_i^\dag-(\Gamma \bar m+\eta)|_Q^2
+|(x_i-x_i^\dag)+\Gamma(\bar m-x^{(N)})|_Q^2\\
&+2 [x_i^\dag -(\Gamma \bar m+\eta)]^TQ [(x_i-x_i^\dag)+\Gamma(\bar m-x^{(N)})]. \end{align*}
 The cost can be rewritten as
\begin{align}
J_i(u_i, \hat u_{-i}, f)=\ &\bar J_i (u_i, f) +{\mathbb E}\int_0^T|(x_i-x_i^\dag)+\Gamma(\bar m-x^{(N)})|_Q^2dt \nonumber\\
&+ {\mathbb E}\left[ | x_i(T)-x_i^\dag (T)|_H^2\right] \nonumber   \\
& +2 {\mathbb  E}\int_0^T\left[x_i^\dag -(\Gamma \bar m+\eta)\right]^TQ \left[(x_i-x_i^\dag)+\Gamma(\bar m-x^{(N)})\right]dt \nonumber \\
&+ 2 {\mathbb  E}\left[ (x_i^\dag(T))^T H (x_i(T)-x_i^\dag(T))\right]  \label{jjb} \\
\le \ & \bar J_i (u_i, f)+C\left(1/N+|x^{(N)}(0)-m_0|^2\right)  \nonumber\\
  &+2 {\mathbb  E}\int_0^T\left[x_i^\dag -(\Gamma \bar m+\eta)\right]^TQ \left[(x_i-x_i^\dag)+\Gamma(\bar m-x^{(N)})\right]dt  \nonumber \\
&+ 2 {\mathbb  E}\left[ (x_i^\dag(T))^T H (x_i(T)-x_i^\dag(T))\right],\label{jNxd}
\end{align}
where the inequality follows from Lemma \ref{lemma:xxbar} and \eqref{xxdag}.
Note that  neither  $x_i-x_i^\dag$ nor $\bar m -x^{(N)}$ in \eqref{jjb} depend on $f$. The terms $x_i^\dag$ and $x_i^\dag-(\Gamma \bar m+\eta)$ are affine in $f$, and $-\bar J_i(u_i, f)$ is convex in $f$ by Lemma \ref{lemma:coer}.  Consequently, it follows from  \eqref{jjb} that  $-J_i(u_i, \hat u_{-i}, f)$ is convex in $f$.   For  $u_i\in {\cal U}_{\hat C_1}$, in analogue to  \eqref{je0}, we obtain
\begin{align}
\bar J_i(u_i, f)\le C-(\epsilon_0/2) \|f\|^2_{L^2} , \label{juf}
\end{align}
where $C$ doest not depend on $u_i$.
We have
\begin{align*}
& \left|{\mathbb  E}\int_0^T\left[x_i^\dag -(\Gamma \bar m+\eta)\right]^TQ \left[(x_i-x_i^\dag)+\Gamma(\bar m-x^{(N)})\right]dt \right|\\
\le\ &\left\{ {\mathbb  E}\int_0^T |x_i^\dag -(\Gamma \bar m+\eta)|^2_Qdt\right\}^{1/2}\cdot \left\{   {\mathbb E}\int_0^T|(x_i-x_i^\dag)+\Gamma(\bar m-x^{(N)})|_Q^2 dt\right\}^{1/2}\\
\le \ & C\left(1/\sqrt{N}+ |x^{(N)}(0)-m_0|\right) (1+\|f\|^2_{L^2})^{1/2}\\
\le \ & C+(\epsilon_0/16) \|f\|^2_{L^2}.
 \end{align*}
 Similarly,
 $$
\left| {\mathbb  E}\left[ (x_i^\dag(T))^T H (x_i(T)-x_i^\dag(T))\right]\right|\le C+(\epsilon_0/16) \|f\|^2_{L^2}.
 $$
 Hence, \eqref{jNxd} gives
\begin{align}
J_i(u_i, \hat u_{-i} , f) \le  C-(\epsilon_0/4) \|f\|^2_{L^2} , \label{jfl2}
\end{align}
where $C$ does not depend on $(N, u_i)$.
 So for  given $u_i\in {\cal U}_{\hat C_1}$,  $J_i(u_i,\hat u_{-i}, f)$ attains a finite supreme at some $f_{u_i}$ since it is a continuous functional of $f$, and by \eqref{jfl2}
 we may further find a constant $\hat C_2$ such that
 \begin{align}
 \sup_{u_i\in {\cal U}_{\hat C_1}  }\|f_{u_i}\|_{L^2} \le \hat C_2. \label{fc2}
 \end{align}
By \eqref{jNxd},
\begin{align}
J_i(u_i, \hat u_{-i}, f)\le\ & \bar J_i (u_i, f) +C(1/N+ |x^{(N)}(0)-m_0|^2)\nonumber
\\
& +C\left(1/N+ |x^{(N)}(0)-m_0|^2\right)^{1/2}  \left({\mathbb  E}
\int_0^T|x_i^\dag -(\Gamma \bar m+\eta)|^2_Qdt\right)^{1/2} \nonumber\\
&+C\left(1/N+ |x^{(N)}(0)-m_0|^2\right)^{1/2}\left( {\mathbb E}
|x_i^\dag(T)|^2\right)^{1/2}. \label{jbjdif}
 \end{align}

Now for $u_i\in  {\cal U}_{\hat C_1}$  and the resulting $f_{u_i}$ satisfying \eqref{fc2}, we further obtain
$$
{\mathbb E} |x_i^\dag(T)|^2 + {\mathbb  E}
\int_0^T|x_i^\dag -(\Gamma \bar m+\eta)|^2_Qdt\le C.
$$
For $u_i\in {\cal U}_{\hat C_1}$,
   \eqref{jbjdif} gives
\begin{align}
\sup_f J_i (u_i, \hat u_{-i}, f) &\le \bar J_i (u_i, f_{u_i})+C(1/\sqrt{N}+ |x^{(N)}(0)-m_0|)\nonumber \\
 & \le \bar J_i (u_i, \hat f_{u_i})    + C(1/\sqrt{N}+ |x^{(N)}(0)-m_0|), \nonumber
\end{align}
where $\hat f_{u_i}$ is determined by Theorem \ref{theorem:p1}. Due to \eqref{juf}, \begin{align}
\sup_{u_i\in {\cal U}_{\hat C_1}}\|\hat f_{u_i}\|_{L^2} \le C  \label{hfC}
\end{align}
 for some constant $C$.
By \eqref{hfC} and the method in \eqref{jjb}, we similarly derive
\begin{align}
 J_i (u_i, \hat u_{-i}, \hat f_{u_i}) \ge \bar J_i (u_i, \hat f_{u_i})-C(1/\sqrt{N}+ |x^{(N)}(0)-m_0|). \nonumber
\end{align}
Hence, for all $ u_i\in {\cal U}_{\hat C_1} $,
$$
\sup_f J_i (u_i, \hat u_{-i}, f)  \ge \bar J_i (u_i, \hat f_{u_i})-C(1/\sqrt{N}+ |x^{(N)}(0)-m_0|).
$$
The constant $C$ in various places does not depend on $u_i$.
The lemma follows. \qed

{\bf  Proof of Theorem \ref{theorem:nash}.} It suffices to show the first inequality by checking $u_i\in {\cal U}_{\hat C_1}$. By Lemma \ref{lemma:jjb}, we have
\begin{align}
\sup_f J_i (u_i, \hat u_{-i}, f) &\ge  \bar J_i(u_i, \hat f_{u_i}) -C_1(1/\sqrt{N}+ |x^{(N)}(0)-m_0|)
\nonumber\\
&\ge \bar J_i(\hat u_i, \hat f_{\hat u_i}) -C_1(1/\sqrt{N}+ |x^{(N)}(0)-m_0|). \label{jlb}
\end{align}
On the other hand, by taking the particular control $\hat u_i$ in Lemma \ref{lemma:jjb},
\begin{align}
 \sup_{f} J_i(\hat u_i, \hat u_{-i}, f)\le \bar J_i(\hat u_i, \hat f_{\hat u_i}) +C_2 (1/\sqrt{N}+ |x^{(N)}(0)-m_0|) . \label{jup}
\end{align}
Subsequently, \eqref{jlb} and \eqref{jup} imply
$$
\sup_f J_i (u_i, \hat u_{-i}, f) \ge \sup_{f} J_i(\hat u_i, \hat u_{-i}, f)
-(C_1+C_2)(1/\sqrt{N}+ |x^{(N)}(0)-m_0|)  .
$$
This completes the proof.  \qed

\section{Further Generalization to Random Initial States}
\label{sec:gen}

This section extends the results to a more general model with random initial states.
For agent ${\cal A}_i$, its dynamics are given by
\begin{equation}
dx_i^o(t)=(Ax_i^o(t)+Bu_i(t)+Gx^{o(N)}(t) +f(t)) dt +D dW_i(t), \quad 1\leq i \leq N, \nonumber 
\end{equation}
where $x^{o(N)}=({1}/{N})\sum_{j=1}^{N}x^o_{j}$.
The initial states of the agents are  given by $x_i^o(0)=\xi_{i}$.
As in \eqref{jcost},  we define $J_i(u_i, u_{-i}, f)$ by using $x_j^o$ in place of $x_j$, $1\le j\le N$. Let $\{{\cal F}_t^o\}_{0\le t\le T}$ be the filtration generated by $\{\xi_i, W_i(t),  1\le i\le N\}$, and $L^2_{{\cal F}^o}(0,T;\mathbb{R}^k)$ is defined accordingly.

({\bf H0}) The sequence $\{\xi_{i}, i\ge 1\}$ consists of independent random variables which are also independent of the Browian motions $\{W_i, i\ge 1\}$. In addition, $\lim_{N\rightarrow \infty}(1/N)\sum_{i=1}^N {\mathbb E}\xi_{i}= m_0$, $\sup_i{\mathbb E}|\xi_{i}|^2\le c_0$ for some constant $c_0$ independent of $N$.

For fixed $\bar u$,
we  consider the  FBSDE
\begin{equation}\label{fbsdeo}
\begin{cases}
dx_i^o =( Ax_i^o +BR^{-1}B^Ty_i^o +G m_i^o+\gamma  p_i^o) dt +D dW_i,\\
 \dot{m_i^o} = (A + G) m_i^o +B\bar u +\gamma  p_i^o,\\
 \dot{ p_i^o} = -(A+G)^T  p_i^o   -(I-\Gamma)^T Q [\mathbb{E}x_i^o- (\Gamma m_i^o+\eta)],\\
dy_i^o=\left\{-A^T y_i^o +Q[x_i^o-(\Gamma m_i^o+\eta)]\right\} dt +\zeta_i^o dW_i,
 \end{cases}
\end{equation}
where $x_i^o(0)=\xi_i$, $m_i^o(0)=m_0$, $ p_i^o(T)=H{\mathbb E} x_i^o(T)
$, and $y_i^o(T)=-Hx_i^o(T)$. Except the random initial state, this FBSDE has the same form as \eqref{fbsdeubar}.

For the current situation where the filtration is not generated only
by the Brownian motions, the proof of Lemma \ref{lemma:uy} is not
applicable.  The solution procedure of (P2) as presented in Section \ref{sec:sub:p2} is only heuristically applied to derive \eqref{fbsdeo}. Nevertheless, we can study  \eqref{fbsdeo} directly and use it to construct decentralized strategies. We still define
$J^{\rm wo}_i(u_i, u_{-i}) =\sup_{f\in L^2(0,T;\mathbb{R}^n)} J_i(u_i, u_{-i} , f)$.
 The next theorem subsumes Corollary \ref{theorem:mpyp2} and  Theorem \ref{theorem:nash}.

\begin{theorem}
Assume that \emph{(H0)}-\emph{(H2)} hold and \eqref{mfge} has a unique solution $({\bf m},{\bf p},{\bf y})$. We further take $\bar u= R^{-1} B^T{\bf y}$ in \eqref{fbsdeo}.
Then the two assertions hold.

\emph{(i)} \eqref{fbsdeo} has  a unique solution in $L^2_{{\cal F}^o}(0,T;\mathbb{R}^n)
\times C^1([0,T];\mathbb{R}^{2n})\times L^2_{{\cal F}^o}(0,T;\mathbb{R}^{2n})$.

\emph{(ii)}
For
$\hat u_i=R^{-1} B^T y_i^o,\ 1\le i\le N,$ we have
\begin{align}
J_i^{\rm wo}(\hat u_i, \hat u_{-i})-\varepsilon_N \le  \inf_{u_i\in {\cal U}} J_i^{\rm wo}(u_i, \hat u_{-i})  \le J_i^{\rm wo}(\hat u_i, \hat u_{-i})
\label{jwo2},
\end{align}
where
 $0\le \varepsilon_N = O(1/\sqrt{N}+ |(1/N)\sum_{j=1}^N {\mathbb E}\xi_j-m_0|)$.
\end{theorem}

\proof (i) Consider \eqref{fbsdeubar} by setting
\begin{equation}
\bar u= R^{-1} B^T {\bf y}, \quad x_i(0)={\mathbb  E}\xi_i.  \label{xoxi}
\end{equation}
Further construct the ODE by taking expectation in \eqref{fbsdeubar}:
\begin{equation}\label{fbsdeoo}
\begin{cases}
\dot{\bar x}_i = A\bar x_i +BR^{-1}B^T\bar y_i +G \bar m_i+\gamma  \bar p_i ,\\
 \dot{\bar m}_i = (A + G) \bar m_i +B\bar u +\gamma  \bar p_i,\\
 \dot{ \bar p}_i = -(A+G)^T  \bar p_i   -(I-\Gamma)^T Q [\bar x_i- (\Gamma \bar m_i+\eta)],\\
\dot{\bar y}_i=-A^T {\bar y}_i +Q[\bar x_i-(\Gamma \bar m_i+\eta)] ,
 \end{cases}
\end{equation}
where $\bar x_i(0)={\mathbb E} \xi_i$, $\bar m_i(0)=m_0$, $ \bar p_i(T)=H\bar  x_i(T)$, and $\bar y_i(T)=-H\bar x_i(T)$. Since \eqref{fbsdeubar} subject to \eqref{xoxi} has a unique solution, \eqref{fbsdeoo} has a solution in $C^1([0,T];\mathbb{R}^{4n})$.
 If \eqref{fbsdeoo} has two different solutions, we will be able to construct two different solutions to \eqref{fbsdeubar} satisfying \eqref{xoxi}, a contradiction to Theorem \ref{theorem:ror}. So \eqref{fbsdeoo} has a unique solution $(\bar x_i,\bar m_i, \bar p_i, \bar y_i)$.

Setting $(m_i^o,p_i^o)= (\bar m_i, \bar p_i)$ in the first and last equations of \eqref{fbsdeo}, we construct the new equations
\begin{equation}\label{fbsdeo2}
\begin{cases}
dx_i^o =( Ax_i^o +BR^{-1}B^Ty_i^o +G \bar m_i+\gamma  \bar p_i) dt +D dW_i,\\
dy_i^o=\left\{-A^T y_i^o +Q[x_i^o-(\Gamma \bar m_i+\eta)]\right\} dt +\zeta_i^o dW_i,
 \end{cases}
\end{equation}
where $x_i^o(0)=\xi_i$ and $y_i^o(T)=-Hx_i^o(T)$.
Let $P$ be  the solution of the Riccati equation \eqref{ricc} and take the transformation $y_i^o= -P x_i^o +\phi$. We obtain
\begin{align}
d\phi = \left[-(A-BR^{-1} B^T P)^T \phi + P(G\bar m_i+\gamma  \bar p_i)-Q(\Gamma \bar m_i+\eta)\right] dt +(\zeta_i^o+PD) dW_i, \nonumber
\end{align}
where $\phi(T)=0$. We solve $(\phi, \zeta_i^o)\in L^2_{\cal F}(0,T; \mathbb{R}^{2n})$, and further obtain $(x_i^o, y_i^o)\in L_{{\cal F}^o}^2(0,T;\mathbb{R}^{2n})$. Subsequently, we can show ${\mathbb E} x_i^o= \bar x_i$. Hence $(x_i^o, m_i^o, p_i^o, y_i^o, \zeta_i^o)$ satisfies \eqref{fbsdeo}. By taking the variation of the first three equations of \eqref{fbsdeo} and applying   an optimal control interpretation as in proving Theorem \ref{theorem:ror}, we can show that
 $(x_i^o, m_i^o, p_i^o, y_i^o, \zeta_i^o)$ is the unique solution.

 (ii) By slightly modifying the proofs of Theorem 18 and the associated lemmas, we can show
 $$
 \sup_{0\le t\le T}{\mathbb E}|\hat{u}^{(N)} -\bar u^*|^2=O(1/N)+O(|{\mathbb E}x^{(N)}(0)-m_0|^2). $$
Next, we adapt the proofs of Lemmas \ref{lemma:hx-bx}-\ref{lemma:jjb} taking into account the random initial states satisfying (H0).
This gives the desired estimate for $\varepsilon_N$. \qed





\section{Conclusion}
\label{sec:con}

This paper introduces a class of mean field LQG games with drift uncertainty. By using the idea of robust optimization, the local strategy is designed by minimizing the worst case cost.   When the decentralized strategies are implemented in a finite population, their performance is characterized as a robust $\varepsilon$-Nash equilibrium.

In this paper we only deal with drift uncertainty. If the Brownian motions are also subject to an uncertain coefficient process to model volatility uncertainty \cite{LPVD83}, the resulting optimal control problems will give a set of more complicated FBSDE. It is also of potential interest to address model uncertainty of the mean field game in a different setup by considering measure uncertainty \cite{CR07,LS07,UP01} in the robust optimization problem. This will necessitate the use of different techniques for analysis. 


\section*{Appendix A}
\renewcommand{\theequation}{A.\arabic{equation}}
\setcounter{equation}{0}
\renewcommand{\thetheorem}{A.\arabic{theorem}}
\setcounter{theorem}{0}

For proving Lemma \ref{lemma:p2cv}, we give another lemma first.
Consider an  auxiliary  optimal control problem with dynamics
 \begin{align}\label{zzqb}
 \left\{\begin{array}{l}
\dot{z_i} = A z_i+ B v_i +G z  +\gamma q,\\
 \dot{z} =(A+G) z +\gamma  q , \\
  \dot{q} = -(A+G)^T q  -(I-\Gamma)^T Q (\mathbb{E}z_i -\Gamma z) ,
\end{array}\right.
\end{align}
where $z_i(0)= z(0)=0$, $q(T)=H{\mathbb E}z_i(T)$ and
$v_i\in L_{\cal F}^2(0,T; \mathbb{R}^{n_1})$. Following the argument in the proof of Lemma \ref{lemma:euzzq}, under (H1) we can show the existence and uniqueness of a solution to \eqref{zzqb}.
The optimal control problem is
\begin{align*}
{\bf (P2b)}  \quad {\rm minimize}\quad \bar J_i^b (v_i) =\mathbb{E}\int_0^T \left\{| z_i - \Gamma z|_Q^2 +v_i^TRv_i -\gamma |q(t)|^2\right\} dt +
 \mathbb{E}| z_i(T)|_H^2.
\end{align*}

Similarly, we may define positive definiteness of $\bar J_i^b$ as in Section \ref{sec:robust}.

\begin{lemma} \label{lemma:jb}
$\bar J_i^a$ is positive semi-definite (resp., positive definite) if and only if $\bar J_i^b$ is positive semi-definite (resp., positive definite).
\end{lemma}

\proof If suffices to show the ``only if" part.

Suppose that $\bar J_i^a$ is positive semi-definite.  Consider any control $v_i\in L_{\cal F}^2(0,T; \mathbb{R}^{n_1})$ for $\bar J_i^b$, and
this gives a unique solution $(z_i, z, q)$.
We take expectation in \eqref{zzqb} to obtain
\begin{equation*}
\begin{cases}
\dot{\bar z}_i = A \bar z_i+ B \bar v_i +G z  +\gamma q,\\
 \dot{z} =(A+G) z +\gamma  q , \\
  \dot{q} = -(A+G)^T q -(I-\Gamma)^T Q (\bar z_i -\Gamma z) ,
\end{cases}
\end{equation*}
where $\bar z_i={\mathbb E }z_i$ and $\bar v_i={\mathbb E} v_i$.

It follows that
\begin{align*}
\bar J_i^b(v_i) &= \bar J_i^a (\bar v_i) +{\mathbb E }\int_0^T \left[|z_i-{\mathbb E}z_i|_Q^2 +| v_i-{\mathbb E}v_i|_R^2\right] dt + {\mathbb E}|z_i(T)-{\mathbb E}z_i(T)|_H^2 \\
 &\ge \bar J_i^a (\bar v_i)\ge 0.
\end{align*}
On the other hand, $\bar J_i^a(0)=0$. This shows that $\bar J_i^b$ is positive semi-definite. The above reasoning is also valid for the positive definite case. This proves the ``only if" part.  \qed

{\it Proof of Lemma \ref{lemma:p2cv}}.
Let $(x_i, m_i, p_i)$  and $( x_i', m_i',  p'_i)$ be two state processes in (P2) corresponding to the controls $u_i$ and $u_i'$, respectively. Assume $\lambda_1\in [0,1]$ and $\lambda_1+\lambda_2=1$.
 We have
\begin{align*}
&\lambda_1\bar J_i (u_i)+\lambda_2 \bar J_i(u_i')- \bar J_i(\lambda_1 u_i+\lambda_2 u_i' )\\
=\ &  \lambda_1\lambda_2 \mathbb{E}\int_0^T \left\{| x_i-x_i' - \Gamma(m_i-m_i')|_Q^2+| u_i-u_i'|_R^2  -\gamma |p_i(t)- p_i'(t)|^2\right\} dt \\
&+\lambda_1\lambda_2 \mathbb{E}| x_i(T)-x_i'(T)|_H^2.
\end{align*}

Denote $z_i=x_i-x_i'$, $z=m_i-m_i'$, $q=p_i- p_i'$ and $v_i= u_i-u_i'$.
It is obvious
$$
\lambda_1\bar J_i (u_i)+\lambda_2 \bar J_i(u_i')- \bar J_i(\lambda_1 u_i+\lambda_2 u_i' )= \lambda_1\lambda_2 \bar J_i^b(v_i).
$$
Recalling Lemma \ref{lemma:jb},
this completes the proof. \qed

\section*{Appendix B}
\renewcommand{\theequation}{B.\arabic{equation}}
\setcounter{equation}{0}
\renewcommand{\thetheorem}{B.\arabic{theorem}}
\setcounter{theorem}{0}

We introduce the FBSDE
\begin{equation} \label{mfge2}
\begin{cases}
dx_i =( Ax_i +BR^{-1}B^Ty_i +G m_i+\gamma  p_i) dt +D dW_i,\\
 \dot{ m_i} = (A + G) m_i +BR^{-1} B^T {\mathbb E} y_i +\gamma  p_i,\\
 \dot{ p_i} = -(A+G)^T p_i   -(I-\Gamma)^T Q [m_i- (\Gamma m_i+\eta)] ,\\
dy_i=\left\{-A^T y_i +Q[x_i-(\Gamma m_i+\eta)]\right\} dt +\zeta_i dW_i,
\end{cases}
\end{equation}
where ${ x}_i(0)=m_i(0)=m_0$, $ p_i(T)=H m_i(T)$, and  $y_i(T)=-Hx_i(T)$. This FBSDE is slightly different from \eqref{mfge} by the third equation and the condition on $ p_i(T)$ and will be more convenient for analysis.

The next lemma shows that the two equation systems  \eqref{mfgeu} and \eqref{mfge2} are equivalent. The proof is straightforward since ${\mathbb E} x_i$ and $m_i$ satisfy the same ODE  with the same initial condition.

\begin{lemma}\label{lemma:sderelate}
If $(x_i,m_i,  p_i, y_i, \zeta_i)\in S[0,T]$ satisfies one of \eqref{mfgeu} and \eqref{mfge2}, it also satisfies the other.~\qed
\end{lemma}

Consider the  ODE system
\begin{equation} \label{mfgode4}
\begin{cases}
\dot{\bar x}_i = A\bar x_i +BR^{-1}B^T\bar y_i +G \bar m_i+\gamma  \bar p_i,\\
 \dot{\bar m}_i = (A + G)\bar m_i +BR^{-1} B^T \bar y_i +\gamma  \bar p_i,\\
 \dot{ \bar p}_i = -(A+G)^T  \bar p_i   -(I-\Gamma)^T Q [\bar m_i- (\Gamma \bar m_i+\eta)] ,\\
\dot{\bar y}_i=-A^T \bar y_i +Q[\bar x_i-(\Gamma \bar m_i+\eta)] ,
\end{cases}
\end{equation}
where $\bar x_i(0)= \bar m_i(0)=m_0$, $\bar p_i(T)=H\bar m_i(T) $ and $\bar y_i(T)= -H\bar x_i(T)$.

\begin{lemma}\label{lemma:sdeode}

The  two statements are equivalent:

\emph{(i)} The \emph{FBSDE}  \eqref{mfge2} has a unique solution in $S[0,T]$.

\emph{(ii)} The \emph{ODE} \eqref{mfgode4} has a unique solution in $C^1([0,T];\mathbb{R}^{4n})$.
\end{lemma}

{\it Proof.} Step 1.   Suppose that (ii) holds and let the unique solution   be denoted by
$(\bar x_i, \bar m_i,  \bar p_i, \bar y_i)$.

Take $(m_i  , p_i )=(\bar m_i, \bar p_i)$ on the right hand side of  the first  and last equations of \eqref{mfge2} to write
\begin{equation} \label{mfg14}
\begin{cases}
dx_i =( Ax_i +BR^{-1}B^Ty_i +G \bar m_i+\gamma \bar p_i) dt +D dW_i,\\
dy_i=\left\{-A^T y_i +Q[x_i-(\Gamma\bar m_i+\eta)]\right\} dt +\zeta_i dW_i,
\end{cases}
\end{equation}
where $y_i(T)=-Hx_i(T)$.
Denote the Riccati equation
\begin{equation}
\dot P+A^TP+PA-PBR^{-1}B^TP+Q=0, \quad P(T)=H,\label{ricc}
\end{equation}
which has a unique solution on $[0,T]$.
Setting  $y_i= -Px_i+\phi$ in \eqref{mfg14}, we obtain two decoupled equations for $(x_i, \phi)$ which is uniquely solved. This further gives a unique solution $(x_i, y_i, \zeta_i)\in L^2_{\cal F}(0,T;\mathbb{R}^{3n})$ for \eqref{mfg14}.
 Taking expectation on both sides of \eqref{mfg14} yields
 \begin{equation} \label{mfg14e}
\begin{cases}
\frac{d}{dt}{\mathbb E}x_i = A{\mathbb E}x_i
 +BR^{-1}B^T{\mathbb E}y_i +G \bar m_i+\gamma \bar p_i ,\\
\frac{d}{dt}{\mathbb E}y_i=-A^T{\mathbb  E}y_i +Q[{\mathbb E}x_i-(\Gamma\bar m_i+\eta)] ,
\end{cases}
\end{equation}
 where ${\mathbb E}y_i(T)=-H{\mathbb E}x_i(T) $.
  By combining \eqref{mfg14e} with the first and fourth equations of \eqref{mfgode4}, it is easy to show  ${\mathbb E}x_i=\bar x_i$ and ${\mathbb E}y_i =\bar y_i$
for all $t\in [0,T]$. This implies
\begin{align*}
\dot{\bar m}_i &= (A+G) \bar m_i +BR^{-1}B^T \bar y_i+\gamma \bar p_i \\
&= (A+G) \bar m_i +BR^{-1}B^T{\mathbb E}  y_i+\gamma \bar p_i.
\end{align*}
The third equation of \eqref{mfge2} is clearly satisfied by $(\bar m_i, \bar p_i)$.
Therefore, $(x_i, m_i,  p_i, y_i,\zeta_i)\defeq (x_i, \bar m_i,  \bar p_i, y_i,\zeta_i)$  satisfies \eqref{mfge2}.

We continue to show that $(x_i,  m_i,  p_i, y_i,\zeta_i)$ above is the unique solution of   \eqref{mfge2}.
Suppose that $(x_i', m_i',  p_i', y_i',\zeta_i')$ is another solution of \eqref{mfge2}.
It is clear that $({\mathbb E}x_i',  m_i',  p_i', {\mathbb E}y_i')$ is a solution of \eqref{mfgode4}. Since \eqref{mfgode4} has a unique solution  $(\bar x_i, \bar m_i,  \bar p_i, \bar y_i)$, we have
$(m_i', p_i')=(\bar m_i,\bar p_i)$. By  using the first and fourth equations of \eqref{mfge2}, we derive the equations satisfied by $(x_i'-x_i,y_i'-y_i )$ and further infer $(x_i', y_i')=(x_i, y_i)$. We conclude that (i) holds.

Step 2.  Suppose that (i) holds with the unique solution denoted by $(x_i, m_i,  p_i, y_i, \zeta_i)$. It is obvious that $(\bar x_i,  \bar m_i, \bar p_i, \bar y_i )\defeq({\mathbb E}x_i, m_i,  p_i, {\mathbb E}y_i)$ is a solution of  \eqref{mfgode4}.  Suppose that $(\bar x_i',  \bar m_i', \bar p_i', \bar y_i')\ne (\bar x_i,  \bar m_i, \bar p_i, \bar y_i )$ is another solution of \eqref{mfgode4}.
Then $(x_i, m_i,  p_i, y_i, \zeta_i)+ (\bar x_i'-\bar x_i,  \bar m_i'-\bar m_i, \bar p_i'-\bar p_i, \bar y_i'-\bar y_i,0)$ is also a solution of \eqref{mfge2}, a contradiction to (i).
So \eqref{mfgode4} has a unique solution. \qed

\begin{lemma} \label{lemma:3to4}
\emph{(i)} If $(\bar x_i, \bar m_i, \bar p_i, \bar y_i )$ is a solution of \eqref{mfgode4},  $ ({\bf m},{\bf p},{\bf y})\defeq (\bar m_i, \bar p_i, \bar y_i)$ satisfies \eqref{mfge}.

\emph{(ii)} If $({\bf  m}, {\bf  p},{\bf  y})$ is a solution of \eqref{mfge}, there exists   $\bar x_i$ such that $(\bar x_i, \bar m_i,  \bar  p_i,  \bar  y_i )\defeq  (\bar x_i, {\bf m},  {\bf p},  {\bf y} )$
satisfies \eqref{mfgode4}.

\emph{(iii)} The \emph{ODE} \eqref{mfgode4} has a unique solution if and only if \eqref{mfge} has a unique solution.
\end{lemma}

\proof (i) If $(\bar x_i, \bar m_i, \bar p_i, \bar y_i )$ is a solution of \eqref{mfgode4}, $\bar x_i=\bar m_i$ and therefore $\bar y_i(T)=-H\bar x_i(T)=-H\bar m_i(T)$. So $({\bf m},{\bf p},{\bf y})$ defined above satisfies \eqref{mfge}.

(ii)  If $({\bf m}, {\bf  p}, {\bf y})$ is a solution of  \eqref{mfge}, we
set $(\bar m_i, \bar p_i,\bar y_i)= ({\bf m},{\bf p},{\bf y})$ and  define $\bar x_i$ by the ODE
$$
\dot{\bar x}_i = A\bar x_i +BR^{-1}B^T\bar y_i +G \bar m_i+\gamma \bar p_i,
$$
where $\bar x_i(0)=m_0$. It can be checked that $\bar m_i=\bar x_i$, which gives $\bar y_i(T) = -H \bar m_i(T)= -H \bar x_i(T)$.
Hence, $(\bar x_i, \bar m_i ,\bar p_i, \bar y_i)$ is a solution to \eqref{mfgode4}.

(iii) Assume that \eqref{mfge} has a unique solution. Let $(\bar x_i, \bar m_i, \bar p_i, \bar y_i)$ and $(\bar x_i', \bar m_i', \bar p_i', \bar y_i') $
be two solutions of \eqref{mfgode4}. By (i), $(\bar m_i, \bar p_i, \bar y_i)$ and $( \bar m_i', \bar p_i', \bar y_i') $ are two solutions of \eqref{mfge} and so must be equal, which further implies $\bar x_i=\bar x_i'$ by the first equation of \eqref{mfgode4}. This shows that \eqref{mfgode4} has a unique solution.

Next assume that \eqref{mfgode4} has a unique solution. Let $({\bf m},{\bf p},{\bf y})$ and $({\bf m}',{\bf p}',{\bf y}')$ be two solutions of \eqref{mfge}.    By (ii), we must have $({\bf m},{\bf p},{\bf y})=({\bf m}',{\bf p}',{\bf y}').$  Therefore, \eqref{mfge} has a unique solution.
 \qed

\end{document}